\documentclass[11pt]{article}
\usepackage{xcolor-material}
\usepackage{amsmath,amsfonts,amscd,bezier}
\usepackage{latexsym,amssymb}
\usepackage[utf8]{inputenc} 
\usepackage[T1]{fontenc}    
\usepackage[english]{babel}
\newtheorem{lemma}{Lemma}[section]
\newtheorem{theorem}[lemma]{Theorem}
\newtheorem{remark}[lemma]{Remark}
\newtheorem{proposition}[lemma]{Proposition}
\newtheorem{corollary}[lemma]{Corollary}

\newtheorem{example}[lemma]{Example}

\newcommand{\cqd}{{\hfill $\rule{2mm}{2mm}$}\vspace{1cm}}
\newcommand{\Dem}{\noindent{\sc Proof:\ \ }}

\begin{document}

	
\title{\bf Bernstein polynomial and \\ value set of differentials for plane branches}
	
\author{{\sc Andrea G. Guimar\~aes and Marcelo E. Hernandes}\thanks{The second-named author was partially supported by CNPq-Brazil.}}
	
\date{}
	
\maketitle
	
\begin{abstract}
In this work, we study the relation between the roots of the Bernstein polynomial $\widetilde{b}$ and the value set of differentials $\Lambda$ for plane branches. For plane branches defined by semiquasihomogeneous polynomial we describe the set of common roots of $\widetilde{b}$ sharing by every branch with a fixed $\Lambda$ set.
\end{abstract}
	
\noindent {\it Keywords: Plane branches, analytic invariants, value set of differentials, Bernstein polynomial.} \\
\noindent {2020 AMS Classification: } 14H20, 32S10, 14F10.\\
	
\section{Introduction}
Let $\mathcal{O}_n=\mathbb{C}\{X\}:=\mathbb{C}\{X_1,\ldots , X_n\}$ be the ring of germs of holomorphic functions at the origin of $\mathbb{C}^n$, and let $\mathcal{D}_X[s]$ denote the ring of differential operators that are holomorphic in $X_1,\ldots ,X_n$ and polynomial in $s$. Given $f\in\mathcal{O}_n$, it is well known (cf. \cite{Bernstein-72} and \cite{Bjork-79}) that
\[I=\{b(s)\in\mathbb{C}[s];\ \exists\ P(s)\in\mathcal{D}_X[s]\ \mbox{such that}\ P(s){f(X)}^{s+1}=b(s){f(X)}^s\}\subset\mathbb{C}[s]\]
is a non-trivial ideal of $\mathbb{C}[s]$. The monic generator $ b_f(s)$ of $I$ is called the {\it Bernstein-Sato polynomial} (or simply, the {\it Bernstein polynomial}) of the germ $f$. Since $b_f(-1)=0$, we consider the reduced Bernstein polynomial $\widetilde{b}_f(s)= b_f(s)/(s+1)$, and denote by $\mathcal{R}_f$ its set of roots, which, according to Kashiwara \cite{Kashiwara-76}, satisfy $\mathcal{R}_f\subset\mathbb{Q}_{<0}$. 

For $f\in\mathcal{O}_n$ with an
isolated singularity at the origin, Malgrange \cite{Malgrange-75} proved that $\sharp\mathcal{R}_f=\mu_f$, where $\mu_f$ is the Milnor number of $f$, and Varchenko \cite{VA} showed that $\mathcal{R}_f$
is an analytic invariant of the hypersurface defined by $f$. 

In what follows, we assume that $f\in\mathcal{O}_n$ defines an
isolated singularity at the origin. 

Associated with the germ $f$, there is another set of $\mu_f$ rational numbers: the set of {\it spectral numbers} $\mathcal{S}_f$. This set is invariant along the stratum $\mu_f$-constant of $f$ (see \cite{VA2}). If $f$ defines an irreducible plane curve (i.e., a plane branch), Saito \cite{Saito-00} provided a complete description of $\mathcal{S}_f$ in terms of the characteristic exponents, or equivalently the minimal generators of the value semigroup $\Gamma_f$ of $f$. 

The relationship between $\mathcal{R}_f=\{\rho_1,\ldots ,\rho_{\mu_f}\}$ and $\mathcal{S}_f=\{\sigma_1,\ldots ,\sigma_{\mu_f}\}$ has been the subject of study by several authors. For instance, in \cite{Ste-76} it was shown that, after a suitable reordering of indices, one has $\sigma_i +\rho_i+1\in \mathbb{N}$ and $\sigma_i+\rho_i\leq n-2$ for all $i=1,\ldots ,\mu_f$. In \cite{andrea}, the first author and A. Hefez proved that $-(\sigma_i +1)\in \mathcal{R}_f$, for all $\sigma_i < \sigma_1+1$, where $\sigma_1:=\min\mathcal{S}_f$ is the smallest spectral number. 

Other authors have also investigated the connection between $\mathcal{R}_f$ and $\mathcal{S}_f$, as well as their relation to the Tjurina number $\tau_f$ of $f$ (e.g. \cite{cassous-2}, \cite{Briancon-89}, \cite{cassous-86}, \cite{cassous-87}, \cite{He-99}, \cite{Saito-89}, \cite{Saito-93}, \cite{Yano}). 

In this paper, we consider a finer analytic invariant for an irreducible plane curve defined by $f\in\mathcal{O}_2$, namely, the set $\Lambda_f\subset\mathbb{N}\setminus\{0\}$ of values of differentials along of $f=0$. 

According to Carbonne \cite{carbonne}, two plane branches $f,h\in\mathcal{O}_2$ are {\it equidifferentiable} if and only if $\Gamma_f=\Gamma_h$ and $\Lambda_f=\Lambda_h$; however, in \cite{osnar} it is shown that $\Lambda_f=\Lambda_h$ already implies $\Gamma_f=\Gamma_h$. Moreover, Zariski \cite{torsion} showed that $\Lambda_f\setminus\Gamma_f=\emptyset$ if and only if $f$ is analytically equivalent to $X_2^{v_0}-X_1^{v_1}$ with $1=gcd(v_0,v_1)<v_0<v_1$, In this case, it is known that $\mathcal{R}_f=\{-(\sigma+1);\ \sigma\in\mathcal{S}_f\}$. The invariant $\Lambda$ plays a central role in the analytic classification of plane branches, as presented in \cite{classification} (or \cite{handbook}). 

More precisely, assuming $\Lambda_f\setminus\Gamma_f\neq\emptyset$, we show (see Theorem \ref{TheoremX}) that for each $\lambda\in\Lambda_f\setminus\Gamma_f$, there exists $\sigma_{\lambda}\in\mathcal{S}_f$, such that $-\sigma_{\lambda}\in\mathcal{R}_f$. In particular, if $\Gamma_f=\langle v_0,v_1\rangle$, i.e., $f$ is analytically equivalent to a semiquasihomogeneous polynomial, we prove (see Theorem \ref{TheoremXX}) that $ \sigma_\lambda = \frac{\lambda}{v_0v_1}$. This gives a proof of the conjecture stated in \cite[Conjecture 1.3]{david}. Futhermore, by fixing an equidifferential class (i.e., a given set $\Lambda_f$), we improve the connection between $\mathcal{S}_f$ and $\mathcal{R}_f$ by proving (see  Corollary \ref{condutor} and Proposition \ref{limitante-inf}) that
$$\begin{array}{ll} 
-(\sigma+1)\in\mathcal{R}_f & \mbox{for every}\ \sigma\in\mathcal{S}_f\ \mbox{such that}\ \sigma<\frac{\lambda_1}{v_0v_1}\ \mbox{where}\ \lambda_1:=\min\Lambda_f\setminus\Gamma_f\ \mbox{and}\\
-\sigma\in\mathcal{R}_f & \mbox{for every}\ \sigma=\frac{\lambda}{v_0v_1}\ \mbox{with}\ \lambda\in\Lambda_f\setminus\Gamma_f.\end{array}$$
In particular, we have $-\sigma\in\mathcal{R}_f$ for every $\frac{\lambda_c}{v_0v_1}\leq\sigma\in\mathcal{S}_f$, where
$\lambda_c:=\min\{\lambda\in\Lambda_f;\ \lambda+\mathbb{N}\subseteq\Lambda_f\}$. This allows us to identify topological types (see Table 1) for which $\Lambda_f$ completely determines $\mathcal{R}_f$. 

The set $\Lambda_f$ is related to the Tjurina number $\tau_f$ by the formula $\tau_f=\mu_f-\sharp\Lambda_f\setminus\Gamma_f$ (see \cite{handbook}). According to Herting and Stahlke (see \cite{He-99}), if we put $\mathcal{E}_f=\{-(\rho+1);\ \rho\in\mathcal{R}_f\}$, that is, the {\it exponents} of $f$, then:
\[\mu_f-\tau_f=\sharp\Lambda_f\setminus\Gamma_f\leq \sum_{\sigma\in\mathcal{S}_f}\sigma-\sum_{\epsilon\in\mathcal{E}_f}\epsilon. \] Note that the difference $\sum_{\sigma\in\mathcal{S}_f}\sigma-\sum_{\epsilon\in\mathcal{E}_f}\epsilon$ counts the number of roots of $\tilde{b}_f$ whose opposites are spectral numbers. This quantity, an analytic invariant, has been studied in several contexts. Bartolo {\it et al.}  \cite{cassous-2} investigated it for plane branches with $\Gamma_f=\langle v_0,v_1,v_2\rangle$. In Example \ref{6-9-v2}, we show that for plane branches with value semigroup $\Gamma_f=\langle 6,9,v_2\rangle$, this quantity exceeds $4$, in line with a comment made by those authors (see page 202 in \cite{cassous-2}).

Since $\mathcal{R}_f$, $\tau_f$ and $\Lambda_f$ are analytic invariants, it is natural to compare the corresponding stratifications in a fixed topological class. Cassous-Nogu\`es (see \cite{cassous-86} and \cite{cassous-87}), Herting and Stahlke (see \cite{He-99}) presented explicit $\tau_f$- and $\widetilde{b}_f$-stratifications for plane branches with value semigroup $\Gamma_f=\langle 6,7\rangle$. In Example \ref{strata}, we enrich this analysis by adding $\Lambda_f$-stratification, thereby illustrating the relationships between the $\tau_f$-, $\widetilde{b}_f$-, and $\Lambda_f$-strata for this topological class.

\section{Brieskorn Lattices and Differential Operators}

In this section we recall some results concerning the satured Brieskorn lattices and identify such a structure with a certain quotient module of differential operators. This module is the main ingredient for proving the Theorem \ref{TheoremX} and Theorem \ref{TheoremXX} which relate the values of differentials along a plane branch with the roots of the reduced Bernstein polynomial $\widetilde{b}_f$.

From now on, we consider a function germ $f\in\mathcal{O}_n$ with an isolated singularity at the origin of $\mathbb{C}^n$, where $n\geq 2$. Let $J(f)=\left(\partial f/\partial X_1,\dots ,
\partial f/\partial X_n\right)$ denote the Jacobian ideal of $f$. The codimension of this ideal as a $\mathbb{C}$-vector space defines the Milnor number $\mu_f:=\dim_{\mathbb{C}}\frac{\mathcal{O}_n}{J(f)}$.

For $m>0$, let $\Omega^m=\Omega^m_{\mathbb{C}^n,0}$ denote the
$\mathcal{O}_n$-module of germs of holomorphic $m$-forms at the
origin of $\mathbb{C}^n$. For $m=0$, we set $\Omega^0=\mathcal{O}_n$.  We consider the following $\mathcal{O}_n$-modules (cf. \cite{Brieskorn-70} or
\cite{Malgrange-74}):
 $$ H'' =\frac{\Omega^n}{df \wedge d\Omega^{n-2}}\ \ \  \ \mbox{and}\ \ \ \ H' = \frac{df \wedge \Omega^{n-1}}{ df \wedge d\Omega^{n-2}},$$
where $df=\sum_{i=1}^{n}f_{X_i}dX_i$. As complex vector spaces, we have
$$ \frac{H''}{H'}\simeq\frac{\Omega^n}{df\wedge \Omega^{n-1}}\simeq \frac{\mathcal{O}_n}{J(f)}.$$ The modules $H''$ and $H'$ can be regarded as $\mathbb{C}\{t\}$-modules via the operation
\begin{equation}\label{t}
t[w] = [fw],\quad \mbox{for}\ [w] \in H''.
\end{equation}

There is a $\mathbb{C}$-linear isomorphism given by
\begin{equation} \label{deltat}
\begin{array}{cccl}
\partial_t: & H' & \longrightarrow & H''\\
                & \left[ df \wedge w\right] &    \mapsto & [dw].\\
\end{array}
\end{equation}
This operator satisfies the Leibniz rule with respect to the
$\mathbb{C}\{t\}$-module structure, and thus uniquely extends to a
meromorphic connection (see \cite{Pham-79}) $$\partial_t :
H''\otimes_{\mathbb{C}\{t\}}\mathbb{C}(t) \longrightarrow
H''\otimes_{\mathbb{C}\{t\}}\mathbb{C}(t).$$
The space $H:=H''\otimes_{\mathbb{C}\{t\}}\mathbb{C}(t)$ is called the
Gauss-Manin system, and $\partial_t$ is the
Gauss-Manin connection associated to $f$ (cf. \cite{Pham-79}).

According to Brieskorn (see \cite{Brieskorn-70}), this connection is regular and
$\text{dim}_{\mathbb{C}(t)}H =\mu_f $. The modules $H'$ and $H''$ are called {\it Brieskorn Lattices}. One also considers the saturation of these lattices with respect to the operator $ \partial_tt$, defined by:
$$\widetilde{H''} = \displaystyle \sum_{k\geq 0} {(\partial_tt)}^kH''\quad \text{and}\quad \widetilde{H'}
= \displaystyle \sum_{k\geq 0} {(\partial_tt)}^kH'.$$ We have that (cf. \cite{Malgrange-75})
\begin{equation}\label{tH2=H1}
   t\widetilde{H''}=\widetilde{H'}
\end{equation}
and the minimal polynomial of
the $ \mathbb{C}$-endomorphism $-\partial_tt$ induced on the quotient $\frac{\widetilde{H''}}{\widetilde{H'}}$, that is,
$$-\partial_tt:
\frac{\widetilde{H''}}{\widetilde{H'}}\rightarrow
\frac{\widetilde{H''}}{\widetilde{H'}},$$
is the reduced
Bernstein polynomial of $f$, denoted by $\tilde{b}_f$. We denote the set of roots of $\widetilde{b}_f$ by \begin{equation}\label{raiz-expoente}\mathcal{R}_f=\{\rho_1,\ldots ,\rho_{\mu_f}\}\ \ \ \mbox{and}\ \ \  \mathcal{E}_f=\{-(\rho+1);\ \rho\in\mathcal{R}_f\}\end{equation} is the exponents of $f$ (that differ by one in the terminology of some authors).

For $ \alpha\in\mathbb{Q}$, we define
$$ \mathcal{V}^{\alpha}H = \bigoplus_{\alpha \leq \beta < \alpha+1 }
\mathbb{C}\{t\} C_{\beta}, \ \ \ \ \
\mathcal{V}^{>\alpha}H=\displaystyle \bigoplus_{\alpha < \beta
< \alpha+1} \mathbb{C}\{t\} C_{\beta},$$ where $
C_{\beta}=\{v\in H;\,\,\, \exists\ l\in\mathbb{N},
{(-\partial_tt+\beta+1)}^l v=0\}$. Thus, the operator $(-\partial_tt+\alpha +1)$ is
nilpotent on the quotient
$\text{Gr}_\mathcal{V}^{\alpha}H:= \mathcal{V}^{\alpha}H/
\mathcal{V}^{>\alpha}H,$ which is isomorphic to $C_{\alpha}$. 

This construction defines the so-called $\mathcal{V}$-{\it filtration}, a decreasing filtration on $H$ such that
$H= \bigcup_{\alpha}\mathcal{V}^{\alpha}H$
and $\mathcal{V}^{\alpha}H $ is a lattice in $H$.

According to Saito (cf. \cite{Saito-89} and \cite{Saito-93}) we get

\begin{proposition}\label{vfiltprop} $ $
\begin{enumerate}
\item \quad $t(\mathcal{V}^{\alpha}H)\subset
\mathcal{V}^{\alpha + 1}H\ \ \mbox{and} \ \
\partial_t(\mathcal{V}^{\alpha}H)\subset \mathcal{V}^{\alpha -
1}H.$
    \item \quad If $v\in \mathcal{V}^\alpha H$, then $hv\in
\mathcal{V}^\alpha H$, for all $h\in \mathcal{O}_n$.
\end{enumerate}
\end{proposition}

If $ N\subseteq M$ are lattices in $H$, then the
$\mathcal{V}$-filtration induces filtrations
on $M$ by 
$$ \mathcal{V}^{\alpha}M=\left(\mathcal{V}^{\alpha}H\right)\cap M \ \ \ \mbox{and}\ \ \ 
\mathcal{V}^{>\alpha}M=\left(\mathcal{V}^{>\alpha}H\right)\cap
M$$ and on the quotient $M/N$ by 
$${\mathcal{V}}^{\alpha}\left(M/N\right)=\left({\mathcal{V}}^{\alpha}M +  N\right)/N\ \ \ \mbox{and}\ \ \
{\mathcal{V}}^{>\alpha}\left(M/N\right)=\left({\mathcal{V}}^{>\alpha}M
+  N\right)/N.
$$
Furthermore:
$$ \text{Gr}_{\mathcal{V}}^\alpha M=\displaystyle\frac{\mathcal{V}^{\alpha} M}{\mathcal{V}^{>\alpha}M}
\simeq \frac{\mathcal{V}^\alpha M+
\mathcal{V}^{>\alpha}H}{\mathcal{V}^{>\alpha}H}\ \ \ \mbox{and}\ \ \ 
\text{Gr}_{\mathcal{V}}^{\alpha}\left(M/N\right)\simeq
\frac{\text{Gr}_{\mathcal{V}}^{\alpha}M}{\text{Gr}_{\mathcal{V}}^{\alpha}N}.$$

Since $H'\subset H''$, we may consider the $\mathcal{V}$-filtration on the quotient $H''/H'$. 
A {\it spectral number} of
$f$ is a rational number $\sigma$ such that
$\text{d}(\sigma) := \text{dim}_{\mathbb{C}} \text{Gr}_{\mathcal{V}}^{\sigma}\left(H''/H'\right)>0.$
We call $\text{d}(\sigma)$ the multiplicity of $\sigma$.
Since $\text{dim}_{\mathbb{C}} \left(H''/H'\right) =
\text{dim}_{\mathbb{C}}\left (\mathcal {O}_n/J(f)\right ) = \mu_f $, 
it follows
that there are exactly $\mu_f$ spectral numbers (counted with its multiplicity). We denote $\mathcal{S}_f$ the set of spectral numbers of $f$. These are invariants under $\mu_f$-constant deformations of $f$ and (cf. \cite[7.3]{SS-85}) satisfy
$$-1<\sigma <n-1\ \ \ \mbox{and}\ \ \ \text{d}(\sigma)=\text{d}(n-2-\sigma)\ \ \ \mbox{for every}\ \ \ \sigma\in\mathcal{S}_f.$$
Moreover, (see \cite[Lemma 3.2.7]{Kulikov-98}) if $\sigma_1:=\min\mathcal{S}_f$ then:
\begin{equation}\label{kulikov}
    H'' \subset\mathcal{V}^{\sigma_1}H''.
\end{equation}

\begin{remark}\label{obs}
We have that $\text{Gr}_{\mathcal{V}}^{\sigma}\left(\widetilde{H''}/\widetilde{H'}\right)\neq
\{0\}$ if and only if $-(\sigma +1)\in\mathcal{R}_f$,
and the multiplicity of this root is the nilpotency degree of
the action of $(-\partial_tt+\sigma + 1)$ on
$\text{Gr}_{\mathcal{V}}^{\sigma}\left(\widetilde{H''}/\widetilde{H'}\right)$
(cf. \cite{Malgrange-75}). Moreover, we get $\mathcal{R}_f\subset\mathbb{Q}\cap (-n,0)$, and if $f\in \mathcal{O}_2$ is irreducible, then $\text{dim}_{\mathbb{C}}\text{Gr}_{\mathcal{V}}^{\sigma}\widetilde{H''}=1$ (see \cite{Saito-89}).
\end{remark}

If $\mathcal{S}_f=\{\sigma_1,\ldots ,\sigma_{\mu_f}\}$ and $\mathcal{E}_f=\{\epsilon_1,\ldots ,\epsilon_{\mu_f}\}$ then after permutation, if necessary, we have (cf. \cite{Ste-76})
\begin{equation}\label{dif-n}
    0\leq \sigma_i-\epsilon_i\leq n-1\ \ \mbox{for}\ \  i=1,\ldots ,\mu_f. 
\end{equation}
Moreover, $\sigma_1=\min\mathcal{S}_f=\min\mathcal{E}_f=\epsilon_1$.

Following the approach of Malgrange in \cite{Malgrange-75} and Brian\c con {\it et al.} in \cite{Briancon-89}, we present some results that play a central role in the proof of Theorem \ref{TheoremX}.

Let $s$ and $U$ be indeterminates, and let $\mathcal{O}_n[f^{-1},s]U$ denote the free $\mathcal{O}_n[f^{-1},s]-$module with basis $\{U\}$. Denoting by $\mathcal{D}_{X,t}$ the ring of differential operators on $\mathcal{O}_n\{t\}$, Malgrange (see \cite{Malgrange-75}) considers  $\mathcal{O}_n[f^{-1},s]U$ as a $\mathcal{D}_{X,t}$-module by defining, for every $g(s)\in \mathcal{O}_n[f^{-1},s]$:
$$
  \partial_{X_i}g(s)U:=\partial_{X_i}(g(s)U):=\frac{\partial g(s)}{\partial X_i}U + sg(s)f^{-1} f_{X_i}U
$$
\begin{equation}\label{acao2}
  tg(s)f^s:= g(s+1)fU
\end{equation}
\begin{equation}\label{acao3}
  \frac{d}{dt}g(s)f^s:=-sg(s-1)f^{-1}U.
\end{equation}

From now on, we replace $U$ by $f^s$ and $f^kU$ by $ f^{s+k}$ for all $k\in \mathbb{Z}$. Consequently, when integer values are assigned to $s$, we recover the usual formulas. On the other hand, if $s$ acts in the standard way, the action of 
$\mathcal{D}$ extends to an action of $\mathcal{D}[s]$ on  $\mathcal{D}_{X,t}f^s$.

It follows by (\ref{acao2}) that
$$t^{-1}g(s)f^s=g(s-1)f^{s-1}.$$ Therefore, multiplication by
$t$ is a bijection on $\mathcal{O}_n[f^{-1},s]f^s$, and
$$-t^{-1}(s+1)g(s)f^s=\frac{d}{dt}g(s)f^s,$$
i.e., multiplication by $s+1$ corresponds to the action of the operator $-t\displaystyle\frac{d}{dt}$. Consequently, using the identity $\displaystyle\frac{d}{dt}t=t\frac{d}{dt}+1$, we see that multiplication by $s$ corresponds to the operator $-\displaystyle\frac{d}{dt}t$.

Given a differential operator $P\in\mathcal{D}_X$ acting on $\mathcal{O}_n$, we consider its representation with coefficients on the right, that is, $P= \displaystyle\sum_{\alpha}\partial_{X}^{\alpha}g_{\alpha}$, where $g_{\alpha}\in \mathcal{O}_n$ and $\partial_{X}^{\alpha}=\partial_{X_1}^{\alpha_1}\dots \partial_{X_n}^{\alpha_n}$ for $\alpha=(\alpha_1,\dots ,\alpha_n)\in\mathbb{N}^n$. 

\begin{proposition}\label{dsfs} Let 
$\mathcal{D}_X[s]\subset \mathcal{D}_{X,t}$ be the polynomial ring in $s=-\frac{d}{dt}t$ with coefficients in $\mathcal{D}_X$. Then:
\begin{displaymath}
{\left({\frac{d}{dt}}\right)}^{-1}\mathcal{D}_X[s]f^s=\left\{\left(Q_0(s)f+\sum_{i=1}^nQ_i(s)f_{X_i}\right)f^s
;\, Q_i(s)\in \mathcal{D}_X[s],\ i=0,\ldots ,n\right\}.
\end{displaymath}
\end{proposition}

\Dem
Let $u\in\mathcal{D}_{X,t}f^s$ such that $\displaystyle
\frac{d}{dt}u\in \mathcal{D}_X[s]f^s$. Then there exists 
$Q(s)\in \mathcal{D}_X[s]$ such that  $\displaystyle\frac{d}{dt}u=
Q(s)f^s$, which implies $(s+1)u=-t\displaystyle
\frac{d}{dt}u=-Q(s+1)f^{s+1}= Q'(s)f^{s+1}$. In particular, for $s=-1$, we obtain
$Q'(-1)(1)=0$. Hence, we must have
$Q'(s)= (s+1)Q_0+\sum_{i=1}^n Q_i\partial_{X_i},$
with $Q_0\in\mathcal{D}_X[s]$ and $Q_i\in\mathcal{D}_X$ for $1\leq i\leq n$. Therefore,
$$(s+1)u= \left((s+1)Q_0+\sum_{i=1}^n Q_i\partial_{X_i}\right)f^{s+1}=(s+1)\left(Q_0f+\sum_{i=1}^n Q_if_{X_i}\right)f^s$$
and consequently
$$u\in \left\{\left(Q_0(s)f+\sum_{i=1}^nQ_i(s)f_{X_i}\right)f^s\, ;\, Q_i(s)\in \mathcal{D}_X[s],\, \, i=0,\dots , n\right\}. $$

To prove the reverse inclusion, consider an element 
$Q(s)f^s =
\left(Q_0(s)f+\sum_{i=1}^nQ_i(s)f_{X_i}\right)f^s\in
\mathcal{D}_X[s]f^s$. Since $ \frac{d}{dt}\partial_{X}=\partial_X\frac{d}{dt}$, it follows from (\ref{acao3}), that
\begin{displaymath}
\begin{array}{ccl}
\displaystyle \frac{d}{dt}Q(s)f^s & = &
-sQ_0(s-1)f^s-\sum_{i=1}^nQ_i(s-1)sf_{X_i}f^{s-1}\\
& = & -sQ_0(s-1)f^s-\sum_{i=1}^nQ_i(s-1)\partial_{X_i}f^s\in
\mathcal{D}_X[s]f^s.
\end{array}
\end{displaymath}
Therefore,
$$\left\{\left(Q_0(s)f+\sum_{i=1}^nQ_i(s)f_{X_i}\right)f^s\, ;\, Q_i(s)\in \mathcal{D}_X[s],\, \,
i=0, \dots , n \right\}\subset \displaystyle
{\left(\frac{d}{dt}\right)}^{-1}\mathcal{D}_X[s]f^s.$$
\cqd

In what follows, we adopt the notations and some results from Section A.1 and Section A.2 of \cite{Briancon-89}. We have chosen to reintroduce some of these results rather than simply citing them, as they require minor corrections and more precise formulations. 

Define 
\begin{displaymath}\label{2.1}
{\xi}_0:={\xi}_0(s)=f^s,
\end{displaymath}
\begin{displaymath}\label{2.2}
\xi_i:=\xi_i(s)=s(s-1)\cdots (s-i+1)f^{s-i},\quad i\in
{\mathbb{N}}\setminus\{0\}.
\end{displaymath}
It then follows from (\ref{acao2}) and 
(\ref{acao3}), that
\begin{displaymath} 
\xi_i= (-1)^i{\left(\frac{d}{dt}\right)}^if^s.
\end{displaymath}

Let $E\subset \mathcal{O}_n$ be a $\mathbb{C}$-vector subspace of dimension $\mu_f$, such that the canonical projection $\mathcal{O}_n\rightarrow \frac{\mathcal{O}_n}{J(f)}$ restricts to a $\mathbb{C}$-linear isomorphism $E\rightarrow \frac{\mathcal{O}_n}{J(f)}$. Then $\mathcal{O}_n=E\oplus J(f)$, and every $g\in \mathcal{O}_n$ can be uniquely written as
$g=e+h,$
where $e\in E$ and $h\in J(f)$. We denote by $DE$ the $\mathbb{C}$-vector subspace generated by ${\partial}^{\alpha}_X e$ with $e\in E$ and $\alpha\in\mathbb{N}^n$.

According to \cite[Proposition A.1.4]{Briancon-89}, if $\mathcal{D}_XJ(f)$ denotes the left ideal in $\mathcal{D}_X$ generated by $J(f)$, then \begin{equation}\label{decomposicao}\mathcal{D}_{X,t}f^s=\mathcal{D}_XJ(f)f^s\oplus \left (\bigoplus_{i\geq 0}DE\xi_i\right ).\end{equation}

\begin{remark}\label{yano-remark} 
Yano (see \cite[Theorem 2.19]{Yano}) shows that the annihilator set {\rm Ann($f^s$)} of $f^s$ in $\mathcal{D}_X$ is precisely the left ideal of $\mathcal{D}_X$ generated by
$${\partial}_{X_i}f_{X_j}-{\partial}_{X_j}f_{X_i},\quad i,j=1,\ldots , n.$$
A direct computation shows that the annihilator ideal {\rm Ann($ \xi_i$)} of $\xi_i$ in $\mathcal{D}_X$ coincides with {\rm Ann($f^s$)}. 
\end{remark}

Now consider $D'$ as the subring of $D$ consisting of all elements whose constant term is zero. We have:
\begin{equation}\label{DE}
     DE\xi_i = D'E\xi_i \oplus E\xi_i.\end{equation}
Indeed, if $ \displaystyle\left(\sum_{\alpha\geq 1} \partial_X^{\alpha}e_{\alpha} + e\right)\xi_i=0$, then, by  Remark \ref{yano-remark}, we conclude that 
$$ \displaystyle \sum_{\alpha  \geq 1} {\partial}_X^{\alpha}e_{\alpha} + e= \sum_{i,j} P_{ij}({\partial}_{X_i}f_{X_j}-{\partial}_{X_j}f_{X_i}), $$
where $P_{ij}\in \mathcal{D}_X$. Rewriting the right-hand side as $\sum_{\alpha}\partial_{X}^{\alpha}g_{\alpha}$, we obtain $e\in J(f)$ and $\sum_{\alpha} \partial_X^{\alpha}e_{\alpha}\in \mathcal{D}_XJ(f)$. Therefore, $\sum_{\alpha} \partial_X^{\alpha}e_{\alpha}=0$ and $e=0$, since $E\cap J(f)=\{0\}$.

Let $\mathcal{I}$ be the $\mathbb{C}$-vector subspace of $J(f)$ generated by 
\begin{equation}\label{bracket} g_{X_l}f_{X_k}- g_{X_k}f_{X_l},\end{equation}
with $g\in \mathcal{O}_n$ and $l,k=1, \ldots , n$. 

Observe that if $\displaystyle \left ( \sum_{|\alpha|\geq 1}{\partial}_{X}^{\alpha}h_{\alpha}+h\right )f^s=0$ with $h_{\alpha}, h\in J(f)$, then $h\in \mathcal{I}$. In fact, by  Remark \ref{yano-remark}, there exist $g_{\gamma_{lk}}, g'_{\gamma_{lk}}\in \mathcal{O}_n$ such that 
\begin{displaymath}
\begin{array}{ccl}
\displaystyle\sum_{|\alpha|\geq 1}{\partial}_{X}^{\alpha}h_{\alpha}+h & = & \sum_{l,k=1}^n\left(\sum_{|\gamma_{lk}|\geq 1}{\partial}_X^{\gamma_{lk}}g_{\gamma_{lk}}\left(  {\partial}_{X_l}f_{X_k}-{\partial}_{X_k}f_{X_l} \right)+ {g'}_{\gamma_{lk}}\left({\partial}_{X_l}f_{X_k}-{\partial}_{X_k}f_{X_l}  \right)\right)\\
& = & \displaystyle\sum_{|\gamma | \geq 1}{\partial}_{X}^{\gamma}g_{\gamma} -{\left(g'_{\gamma_{lk}}\right)}_{X_l}f_{X_k} + {\left(g'_{\gamma_{lk}}\right)}_{X_k}f_{X_l}.  
\end{array}
\end{displaymath}   
By the uniqueness of the representation of the elements in $\mathcal{D}_X$, we obtain 
$$h= {\left(g'_{\gamma_{lk}}\right)}_{X_k}f_{X_l}-{\left(g'_{\gamma_{lk}}\right)}_{X_l}f_{X_k}\in\mathcal{I}.$$

This property, together with (\ref{DE}), ensures the well-definedness of the surjective $\mathbb{C}$-linear map 
\begin{equation} \label{c}
c:\mathcal{D}_{X,t}f^s=\mathcal{D}_XJ(f)f^s\oplus \left (\bigoplus_{i\geq
0}DE\xi_i\right )\,\, \longrightarrow \,\, \displaystyle\frac{J(f)f^s}{\mathcal{I}f^s} \oplus
\left (\bigoplus_{i\geq 0}E\xi_i\right ),
\end{equation}
given by
\begin{center}
$c\left(\left(\sum_{|\alpha|\geq 1}{\partial}_X^\alpha g_{\alpha}+h\right)f^s\right)= \overline{hf^s}\in \displaystyle\frac{J(f)f^s}{\mathcal{I}f^s}$, with $g_{\alpha},\, h\in J(f)$;
\end{center}
and
\begin{center}
$c\left(\left(\sum_{|\alpha|\geq 1}{\partial}_X^\alpha e_{\alpha}+ e\right)\xi_i\right)= e\xi_i\in E\xi_i$, with  $e_{\alpha},\, e\in E$.
\end{center}

\begin{proposition}\label{kernel} Let $c$ be the $\mathbb{C}$-linear map defined in (\ref{c}). Then:
$$\text{{\rm ker}}\,c = \sum {\partial}_{X_i}\left( \mathcal{D}_XJ(f)f^s\oplus \left (\bigoplus_{i\geq 0}DE\xi_i\right )\right).$$
\end{proposition}

\Dem It is immediate that $ \sum {\partial}_{X_i}\left( \mathcal{D}_XJ(f)f^s\oplus \left (\bigoplus_{i\geq 0}DE\xi_i\right )\right)\subseteq \text{ker}\,c$. 

Conversely, suppose 
$$ \displaystyle\Big(\sum_{|\alpha|\geq 1}{\partial}_X^{\alpha}g_{\alpha}+h\Big)f^s + \sum_i\Big(\sum_{|\gamma_i|\geq 1}{\partial}_{X}^{\gamma_i}e_{\gamma_i} + e_i\Big)\xi_i\in \text{ker}\,c.$$ Then $ hf^s \in \mathcal{I}f^s$ and $\displaystyle \sum_ie_i\xi_i=0$. Thus, $e_i\xi_i=0$ and, by  Remark \ref{yano-remark}, we have $e_i=0$ for all $i$. Futhermore, since $\mathcal{I}$ is generated by expressions of the form (\ref{bracket}), it follows that
$$hf^s= \displaystyle\sum_{l,k=1}^n ({(g_{lk})}_{X_l}f_{X_k}- {(g_{lk})}_{X_k}f_{X_l})f^s=\displaystyle\sum_{l,k=1}^n \left(\partial_{X_l}g_{lk}f_{X_k}- \partial_{X_k}g_{lk}f_{X_l}\right)f^s.$$
Hence, $\text{ker}\,c \subseteq 
\sum {\partial}_{X_i}\left( \mathcal{D}_XJ(f)f^s\oplus (\bigoplus_{i\geq 0}DE\xi_i)\right)$.
\cqd

It follows, by Proposition \ref{dsfs}, that
$\mathcal{D}_X J(f)f^s\subseteq
\displaystyle{\left(\frac{d}{dt}\right)}^{-1}\mathcal{D}_X[s]f^s
\subseteq \mathcal{D}_X[s]f^s\ (\subseteq \mathcal{D}_{X,t}f^s)$. Thus, we obtain $\mathbb{C}$-vector spaces $A_2\subseteq A_1\subseteq \left (\bigoplus_{i\geq 0}E\xi_i\right )$ such that
\begin{equation}\label{conjA}\begin{array}{cc}\frac{J(f)f^s}{\mathcal{I}f^s} \oplus A_1=c(\mathcal{D}_X[s]f^s)\ \ \mbox{and} \vspace{0.3cm}\\
\frac{J(f)f^s}{\mathcal{I}f^s} \oplus A_2= c\left({\left(\frac{d}{dt}\right)}^{-1}\mathcal{D}_X[s]f^s\right).\end{array}\end{equation}
Moreover, we have:
$$\displaystyle\frac{J(f)f^s}{\mathcal{I}f^s} \oplus A_2=\{c(v);\,v\in \mathcal{D}_X[s]f^s\ \text{and}\
(s+1)v\in \mathcal{D}_X[s]f^{s+1}\}.$$

Considering the restriction of $c$ to $\mathcal{D}_X[s]f^s$, we obtain the epimorphism
\begin{displaymath}
\begin{array}{cccc}
c_1:& \mathcal{D}_X[s]f^s & \rightarrow & \displaystyle\frac{J(f)f^s}{\mathcal{I}f^s}\oplus A_1 \\
& v& \mapsto & c(v).
\end{array}
\end{displaymath} 
Since, by Proposition \ref{kernel}, we have $\text{ker}\, c = \sum {\partial}_{X_i}\left( \mathcal{D}_XJ(f)f^s\oplus (\bigoplus_{i\geq 0}DE\xi_i)\right)$, it follows that
$$\text{ker}\, c_1 = \Big(\sum {\partial}_{X_i}\big( \mathcal{D}_XJ(f)f^s\oplus (\bigoplus_{i\geq 0}DE\xi_i)\big)\Big)\cap \mathcal{D}_X[s]f^s= \sum {\partial}_{X_i}\mathcal{D}_X[s]f^s .$$
Thus, we obtain the $ \mathbb{C}$-isomorphism 
$$
\begin{array}{cccc}
 \theta_1 : \displaystyle \frac{\mathcal{D}_X[s]f^s}{\sum \partial_{X_i}\mathcal{D}_X[s]f^s} & \longrightarrow & \displaystyle\frac{J(f)f^s}{\mathcal{I}f^s} \oplus A_1  & =c(\mathcal{D}_X[s]f^s)\vspace{0.2cm}\\
 \left[g\right] & \longmapsto & c(g). &
\end{array}
$$

Similarly, restricting $c$ to
${\left (\frac{d}{dt}\right)}^{-1}\mathcal{D}_X[s]f^s$, we obtain the epimorphism
\begin{displaymath}
\begin{array}{cccc}
c_2:& {\left ( \frac{d}{dt}\right)}^{-1}\mathcal{D}_X[s]f^s & \rightarrow & \displaystyle\frac{J(f)f^s}{\mathcal{I}f^s}\oplus A_2 \\
& v& \mapsto & c(v)
\end{array}
\end{displaymath} 
with
$$\text{ker}\, c_2= \Big(\sum {\partial}_{X_i}\big( \mathcal{D}_XJ(f)f^s\oplus (\bigoplus_{i\geq 0}DE\xi_i)\big)\Big)\cap {\left (\frac{d}{dt}\right)}^{-1}\mathcal{D}_X[s]f^s= \sum {\partial}_{X_i}{\left ( \frac{d}{dt}\right)}^{-1}\mathcal{D}_X[s]f^s .$$ This yields the $\mathbb{C}$-isomorphism 
\begin{equation}\label{isodsmodulo11}
\begin{array}{cccc}
 \theta_2 : \displaystyle \frac{{\left ( \frac{d}{dt}\right)}^{-1}\mathcal{D}_X[s]f^s}{\sum \partial_{X_i}{\left ( \frac{d}{dt}\right)}^{-1}\mathcal{D}_X[s]f^s} & \longrightarrow & \displaystyle\frac{J(f)f^s}{\mathcal{I}f^s} \oplus A_2  & =c\left ({\left ( \frac{d}{dt}\right)}^{-1}\mathcal{D}_X[s]f^s  \right )\vspace{0.2cm}\\
 \left[g\right] & \longmapsto & c(g). &
\end{array}
\end{equation}

Now we present some results that allow us to identify the satured Brieskorn lattices in terms of differential operator modules.

\begin{proposition}\label{lambda1}{(\cite[Proposition 5.3]{Malgrange-75})} Let $\Im$ denote the left annihilator of $f^s$ in $\mathcal{D}_X[s]$. Then the epimorphism
$$\begin{array}{cccc}
\psi_1:& \Omega^n[s] & \longrightarrow & \widetilde{H''}\\
          & \sum_{i=0}^r s^ig_idX & \longmapsto & \sum_{i=0}^r {(-\partial_tt)}^i[g_idX]
\end{array}$$ 
induces a $\mathbb{C}$-isomorphism 
$\displaystyle \frac{\Omega^n[s]}{\Omega^n[s]\Im }\longrightarrow\widetilde{H''}$, where 
$$ \Omega^n[s]\Im = \{ \omega(s)P(s);\,\,\, \omega(s)\in \Omega^n[s]\,\,\,\text{and}\,\,\, P(s)\in \Im \}.$$
\end{proposition}

If $P\in\mathcal{D}_X$, the {\it adjoint operator} $P^*$ is defined by the following rules:
\begin{itemize}
    \item[i)] $P^*:=P$, if $P\in\mathcal{O}_n$;
    \item[ii)] $P^*:=-P$, if $P=\partial_{X_i}$;
    \item[iii)] $P^*:=P_2^*P_1^*$, if $P=P_1P_2$ with $P_1, P_2\in \mathcal{D}_X$.
\end{itemize}

Given
$$P(s)=\sum_{i=0}^{l}P_is^i\in \mathcal{D}_X[s]\ \ \ \mbox{and}\ \ \ \omega(s)=\sum_{j=0}^{r}g_jdX\in \Omega^n[s],$$
where $P_i\in \mathcal{D}_X$, $g_i\in\mathcal{O}_n$ and $dX=dX_1\wedge \dots \wedge dX_n$, we define:
$$ \omega(s)P(s):= \displaystyle\sum_{{i=0,\ldots , l}\atop {j=0,\ldots , r}}s^{i+j}P_i^*(g_j)dX.$$ 

\begin{proposition}\label{lambda2}
With the above notation, the map 
\begin{displaymath}
\begin{array}{cccl}
\psi_2:& \displaystyle \frac{\mathcal{D}_X[s]f^s}{\sum \partial_{X_i}\mathcal{D}_X[s]f^s} &  \longrightarrow & \frac{\Omega^n[s]}{\Omega^n[s]\Im }\vspace{0.3cm}\\
      &\overline{(\sum_{i=0}^rP_is^i)f^s} & \longmapsto &  \overline{dX(\sum_{i=0}^rP_is^i)}
\end{array}
\end{displaymath}
is a $\mathbb{C}$-isomorphism.
\end{proposition}

\Dem Note that if $P(s)f^s=Q(s)f^s$ then $P(s)-Q(s)\in \Im $. Therefore, the map 
\begin{displaymath}
\begin{array}{cccl}
\Pi : & \mathcal{D}_X[s]f^s & \longrightarrow & \displaystyle \frac{\Omega^n[s]}{\Omega^n[s]\Im }\vspace{0.2cm}\\
 & (\sum_{i=0}^rP_is^i)f^s & \longmapsto & \overline{dX(\sum_{i=0}^rP_is^i)}
\end{array}
\end{displaymath}
is a surjective $\mathbb{C}$-linear map, and we have $\sum \partial_{X_i}\mathcal{D}_X[s]f^s\subseteq \text{ker}\,\Pi $. 

Conversely, suppose $P(s)f^s\in \text{ker}\, \Pi $. Then there exist $g_i=g_{i0}+sg_{i1}+\cdots + s^{l_i}g_{il_i} \in \mathcal{O}_n[s]$ and $Q_i(s)= Q_{i0} + Q_{i1}s+ \cdots +
Q_{ik_i}s^{k_i}\in \Im $ such that  
$$dXP(s)=\displaystyle\sum_i g_i(s)dXQ_i(s)=dX\sum_ig_i(s)Q_i(s)=dX\displaystyle\sum_i\Big(\sum_{{j=0,\ldots , l_i}\atop{q=0,\ldots , k_i}}g_{ij}Q_{iq}s^{j+q}\Big).$$
Hence,
$$\Big(\sum_{i=0}^rs^iP_i^* - \displaystyle\sum_i \Big(\sum_{{j=0,\ldots , l_i}\atop{q=0,\ldots , k_i}}s^{j+q}{(g_{ij}Q_{iq})}^*\Big)\Big)(1)dX =0$$
so that $$P(s)-\displaystyle\sum_i g_i(s)Q_i(s)=\displaystyle\sum_i \partial_{X_i}R(s),$$
for some $R(s)\in \mathcal{D}_X[s]$. Since $\sum_ig_i(s)Q_i(s)f^s=0$, it follows that 
$$P(s)f^s= \sum_i \partial_{X_i}R(s)f^s.$$

Therefore, $\text{ker}\, \Pi= \sum \partial_{X_i}\mathcal{D}_X[s]f^s$, and $\psi_2$ is the induced isomorphism by $\Pi$.
\cqd 

As an immediate consequence of Proposition \ref{lambda1} and Proposition \ref{lambda2}, we obtain the following result.

\begin{proposition}\label{isomorfismo}
The map
\begin{equation}\label{isdmo}
\vartheta:  \displaystyle \frac{\mathcal{D}_X[s]f^s}{\sum \partial_{X_i}\mathcal{D}_X[s]f^s}   \longrightarrow  \widetilde{H''}
\end{equation}
defined by 
$$ \vartheta\left(\overline{(\sum_{i=0}^rP_is^i)f^s}\right) = \sum_{i=0}^r{(-\partial_tt)}^i[P_i^*(1)dX],$$ 
is a $\mathbb{C}$-isomorphism.
\end{proposition}

Considering the $\mathbb{C}$-isomorphism in (\ref{isdmo}), we obtain the following corollary:

\begin{corollary}\label{isoHlinhatil}
By restriction, the map $\vartheta$ induces the $\mathbb{C}$-isomorphism
$$\displaystyle\frac{
{\left(\frac{d}{dt}\right)}^{-1}\mathcal{D}_X[s]f^s}{\sum
\partial_{X_i}{\left(\frac{d}{dt}\right)}^{-1}\mathcal{D}_X[s]f^s}\cong
\widetilde{H'}.$$
\end{corollary}

\Dem Let $$\left(\sum_{i=0}^rP_ifs^i +
\sum_{i=1}^n\sum_{j=0}^l Q_{ij}f_{X_i}s^j\right)f^s \in
\displaystyle
{\left(\frac{d}{dt}\right)}^{-1}\mathcal{D}_X[s]f^s\subseteq
\mathcal{D}_X[s]f^s$$ with $P_i,\, Q_{ij} \in \mathcal{D}_X$. Then
$$\hspace{-9.5cm} \vartheta \left(\overline{\Big(\sum_{i=0}^rP_ifs^i +  \sum_{i=1}^n\sum_{j=0}^l Q_{ij}f_{X_i}s^j\Big)f^s}\right) = $$
\begin{equation}\label{somatorio}
\begin{array}{l} = \sum_{i=0}^r{(-\partial_tt)}^i[fP_i^*(1)dX]+ \sum_{j=0}^l {(-\partial_tt)}^j\left[\sum_{i=1}^nf_{X_i}Q_{ij}^*(1)dX\right] \\
 =\sum_{i=0}^rt{(-\partial_tt-1)}^i[P_i^*(1)dX]+\sum_{j=0}^l {(-\partial_tt)}^j\left[\sum_{i=1}^n df \wedge {(-1)}^{i+1}Q_{ij}^*(1)dX_1
\wedge \cdots \wedge \widehat{dX_i}\wedge \cdots \wedge dX_n
\right]\end{array}\end{equation}
$$\hspace{-13.5cm}\in t\widetilde{H''} + \widetilde{H'}= \widetilde{H'},$$
where the last equality follows from $t\widetilde{H''}= \widetilde{H'}$ (see (\ref{tH2=H1})). Therefore, $\vartheta \left(\overline{{\left(\frac{d}{dt}\right)}^{-1}\mathcal{D}_X[s]f^s}\right)\subseteq \widetilde{H'}$.
Since any element of $\widetilde{H'}$ can be expressed as a sum as (\ref{somatorio}), we conclude that 
$\vartheta \left(\overline{{\left(\frac{d}{dt}\right)}^{-1}\mathcal{D}_X[s]f^s}\right)= \widetilde{H'}.$
\cqd

In this way, we obtain the following theorem.

\begin{theorem}\label{TheoremX} Given $\overline{0}\neq\overline{g}\in\frac{(f)+J(f)}{J(f)}$, there exists $\sigma \in \mathbb{Q}$, with $\sigma \geq \sigma_1+1$, such that $\overline{0}
\neq \overline{[gdX]}\in \text{Gr}_{\mathcal{V}}^{\sigma}\left(H''/H'\right)$ for some $\sigma \in \mathcal{S}_f$. If $f\in\mathcal{O}_2$ is irreducible, then $-\sigma \in \mathcal{R}_f$. 
\end{theorem}

\Dem Given $gdX=(Q_0f+\sum_{i=1}^{n}Q_if_{X_i})dX$, where $g\notin J(f)$ and $Q_i\in\mathcal{O}_n$ we consider the element $\eta=\sum_{i=1}^{n}(-1)^{i-1}Q_idX_1\wedge\ldots\wedge \widehat{dX_i}\wedge\ldots \wedge dX_n\in\Omega^{n-1}$. Then we get $d\eta=\left (\sum_{i=1}^{n}(Q_i)_{X_i}\right )dX$, and it follows from (\ref{deltat}) that
$$\left[\left (\sum_{i=1}^{n}Q_if_{X_i}\right )dX\right]=[df\wedge \eta] =\partial_t^{-1}[d\eta]=\partial_t^{-1}\left[\left (\sum_{i=1}^{n}(Q_i)_{X_i}\right )dX\right].$$
In this way, by (\ref{t}), we obtain 
$$[gdX]=\left[(Q_0f+\sum_{i=1}^{n}Q_if_{X_i})dX\right]= t\left[Q_0dX\right]+\partial_t^{-1}\left[\left (\sum_{i=1}^{n}(Q_i)_{X_i}\right )dX\right].$$   

Since $g\notin J(f)$, we have $gdX\notin df\wedge \Omega^{n-1}, df \wedge d\Omega^{n-2}$. From (\ref{kulikov}) we know that $[dX]\in\mathcal{V}^{\sigma_1}H''$, thus, by Proposition \ref{vfiltprop}, it follows that $[gdX]\in \mathcal{V}^{\sigma_1+1}H''$. Hence, $\overline{0}\neq \overline{[gdX]}\in \text{Gr}_{\mathcal{V}}^{\sigma}\left(H''/H'\right)$, for some $\sigma \geq \sigma_1+1$ with $\sigma\in\mathcal{S}_f$.

Now consider $f\in\mathcal{O}_2$ irreducible. By (\ref{raiz-expoente}) and (\ref{dif-n}), we have that either $-(\sigma+1)$ or $-\sigma$ belongs to $\mathcal{R}_f$. 

Suppose that $-(\sigma+1)\in\mathcal{R}_f$. Then by Remark \ref{obs}, it follows that $\{0\}\neq \text{Gr}_{\mathcal{V}}^{\sigma}\left (\frac{\widetilde{H''}}{\widetilde{H'}}\right )\simeq\frac{\text{Gr}_{\mathcal{V}}^{\sigma}\widetilde{H''}}{\text{Gr}_{\mathcal{V}}^{\sigma}\widetilde{H'}}$and that $\dim_{\mathbb{C}}\text{Gr}_{\mathcal{V}}^{\sigma}\widetilde{H''}=1$. Since $[gdX]\in H''\subseteq\widetilde{H''}$, we must have $[gdX]\not\in\widetilde{H'}$. However, by Proposition \ref{dsfs},  we have $gf^s\in {\left(\frac{d}{dt}\right)}^{-1}\mathcal{D}_X[s]f^s$ and it follows from (\ref{conjA}) that 
$$
c_2(gf^s)=c(gf^s)\in c\left({\left(\frac{d}{dt}\right)}^{-1}\mathcal{D}_X[s]f^s\right)=\displaystyle\frac{J(f)f^s}{\mathcal{I}f^s}\oplus A_2.$$
Then, by (\ref{isodsmodulo11}) and Corollary \ref{isoHlinhatil}, we have 
$$ \vartheta\circ\theta_2^{-1}(c(gf^s))=[gdX]\in \widetilde{H'}.$$
This yields a contradiction. Hence $-(\sigma+1)\not\in\mathcal{R}_f$, and consequently we must have $-\sigma\in\mathcal{R}_f$.
\cqd

\section{Plane branches}\label{sec:plana}

In Theorem \ref{TheoremX}, we show that if $f\in\mathcal{O}_2$ is irreducible, then for every nonzero element in $\frac{(f)+J(f)}{J(f)}$, there exists $\sigma\in\mathcal{S}_f$ such that $-\sigma\in\mathcal{R}_f$, or equivalently, $\sigma-1\in\mathcal{E}_f$ (see (\ref{raiz-expoente})). Moreover, for any hypersurface with an isolated singularity, Hertling and Stahlke show (cf. \cite[Proposition 3.5]{He-99}) that  
\begin{equation}\label{mu-tau}
\dim_{\mathbb{C}}\frac{\widetilde{H''}}{H''}=\sum_{\sigma\in\mathcal{S}_f}\sigma-\sum_{\epsilon\in\mathcal{E}_f}\epsilon\geq\dim_{\mathbb{C}}\frac{(f)+J(f)}{J(f)}=\mu_f-\tau_f,\end{equation}
where $\tau_f=\dim_{\mathbb{C}}\frac{\mathcal{O}_2}{(f)+J(f)}$ is the Tjurina number of $f$. In particular, $\dim_{\mathbb{C}}\frac{\widetilde{H''}}{H''}$ gives the number of roots of $\tilde{b}_f$ whose opposite are spectral numbers.

From now on, we consider $n=2$, that is, $\mathcal{O}_2=\mathbb{C}\{X_1,X_2\}$. In addition, we assume that $f\in\mathcal{O}_2$ is irreducible with multiplicity $mult(f)>1$; in this way, $f$ defines a germ of an irreducible plane curve singularity $\mathcal{C}_f$, also known as a {\it plane branch}.

Up to analytical equivalence, a plane branch $\mathcal{C}_f$ 
admits a primitive parameterization\footnote{A parameterization is primitive if the exponents have no common divisor greater than $1$.} 
\begin{equation}\label{param}
\varphi(T)=(\varphi_1(T),\varphi_2(T))=\left ( T^{v_0},T^{v_1}+\sum_{i>v_1}a_iT^i\right )\in\mathbb{C}\{T\}\times\mathbb{C}\{T\}\end{equation}
where $v_1>v_0=mult(f)$ and $v_0\nmid v_1$. 

If $\mathcal{O}_f:=\frac{\mathcal{O}_2}{(f)}=\mathbb{C}\{x_1,x_2\}$ denotes the local ring of $\mathcal{C}_f$, then the {\it value semigroup} of $\mathcal{C}_f$ is defined by
\begin{equation}
    \label{def-semi}
\Gamma_f=\left \{\nu_f(h):=ord_T(h(\varphi(T)));\ h\in\mathcal{O}_f\setminus \{0\}\right \}=\{\nu_f(q)=ord_T(q(\varphi(T)));\ q\in\mathcal{O}_2\setminus (f)\}\subset\mathbb{N}\end{equation} and it characterizes the topological type of the plane branch (see \cite{zariski}).

We denote by $\{v_0, v_1,\ldots ,v_g\}$ the minimal set of generators of $\Gamma_f$, and write $\Gamma_f=\langle v_0, v_1,\ldots ,v_g\rangle$, where $v_0<v_1<\ldots <v_g$ and $\gcd(v_0,v_1,\ldots ,v_g)=1$.  
We can compute the value semigroup $\Gamma_f$ from a parameterization of $\mathcal{C}_f$, as (\ref{param}), in the following way:  

Define\begin{equation}\label{expoentes-caract}\begin{array}{ll}\beta_0:=v_0, & e_0:=v_0;\\ \beta_1:=v_1, & e_1=gcd(e_0,\beta_1);\\ \beta_j:=\min\{i;\ a_i\neq 0\ \mbox{and}\ e_{j-1}\nmid i\}, & e_j:=gcd(e_{j-1},\beta_j)\ \ \ \mbox{for}\ \ j>1. \end{array} \end{equation}

Since the parameterization (\ref{param}) is primitive, there is an integer $g'\geq 1$ such that $e_{g'}=1$. Moreover, we have $g'=g$. The integers $\beta_0<\beta_1<\cdots <\beta_g$ are the {\it characteristic exponents} of $\mathcal{C}_f$, and they are related with the generators of $\Gamma_f$ by the relations 
\begin{equation}\label{exp-ger}v_0=\beta_0,\ \ \ v_1=\beta_1,\ \ \ v_{i+1}=n_iv_i+\beta_{i+1}-\beta_i\ \ \ \mbox{with}\ \ \ n_i=\frac{e_{i-1}}{e_i}\ \ \ \mbox{for}\ \ \ 0<i<g.\end{equation}

The semigroup $\Gamma_f$ admits a conductor $c_f:=\min\{\gamma\in\Gamma_f;\ \gamma+\mathbb{N}\subseteq\Gamma_f\}$ which coincides with the Milnor number $\mu_f$ of $f$. According to Zariski (cf. \cite{zariski}), we have $\mu_f=c_f=\sum_{i=1}^{g}(n_i-1)v_i-v_0+1$.

\begin{remark}\label{escrita}
Let $\Gamma_f=\langle v_0,\ldots ,v_g\rangle$ be the value semigroup of a plane branch. Any $z\in\mathbb{Z}$ can be written uniquely as $z=\sum_{i=0}^{g}s_iv_i$ with $s_i\in\mathbb{Z}$ and $0\leq s_i<n_i$ for $1\leq i\leq g$. Moreover, for such an expression, $z\in\Gamma_f$ if and only if $s_0\geq 0$.	\end{remark}

Using the value semigroup, we can obtain a topological lower bound for $\mu_f-\tau_f$ as follows.

\begin{proposition}\label{bound}
Let $f\in\mathcal{O}_2$ be an irreducible element with value semigroup $\Gamma_f=\langle v_0,v_1,\ldots ,v_g\rangle$. Then
\[\dim_{\mathbb{C}}\frac{\widetilde{H''}}{H''}\geq \mu_f-\tau_f\geq \frac{3n_g-2}{4}\left ( \sum_{i=1}^{g-1}(n_i-1)\frac{v_i}{n_g}-\frac{v_0}{n_g}+1\right ). \]
\end{proposition}
\Dem By Corollary 4.10 in \cite{osnar2}, we have 
$\tau_f\leq \mu_f-\frac{3n_g-2}{4}\mu_{g-1}$, where $\mu_{g-1}$ is the condutor of the semigroup $\langle \frac{v_0}{n_g},\ldots ,\frac{v_{g-1}}{n_g}\rangle$, that is, $\mu_{g-1}=\sum_{i=1}^{g-1}(n_i-1)\frac{v_i}{n_g}-\frac{v_0}{n_g}+1$. Thus, the result follows from (\ref{mu-tau}).
\cqd

For $g=2$, Bartolo {\it et al.} in \cite[Proposition 3.2]{cassous-2} show that $\dim_{\mathbb{C}}\frac{\widetilde{H''}}{H''}\geq (n_2-1)\left ( \frac{v_0}{n_2}-1 \right )\left (\frac{v_1}{n_2}-1 \right )$. They show that equality holds for any plane branch with value semigroup $\langle 4,6,v_2\rangle$ and they exhibit a plane branch with value semigroup $\langle 8,10,47\rangle$ achieving the equality. Notice that in both cases $n_2=2$, and the lower bounds for  $\dim_{\mathbb{C}}\frac{\widetilde{H''}}{H''}$  given in \cite{cassous-2} and in Proposition \ref{bound} coincide. 

Considering the topological data, the above lower bounds allow us to conclude that for any plane branch with value semigroup $\Gamma =\langle 6,9,v_2\rangle$, we have $\dim\frac{\widetilde{H''}}{H''}\geq 4$. However, according to \cite[page 202]{cassous-2}, it was not known whether equality is achieved. In fact, this is not the case, as we will show in Example \ref{6-9-v2}.   

The difference $\mu_f-\tau_f$ is related to an analytic invariant of the branch defined by $f$, given by values of differentials as we now explain.

The K\"ahler differential module of $\mathcal{C}_f$ (or $\mathcal{O}_f$) is the $\mathcal{O}_f$-module $\Omega_f:=\frac{\Omega^1}{f\Omega^1+\mathcal{O}_2df}$. If $\varphi(T)=(\varphi_1(T),\varphi_2(T))$ is a parameterization of $\mathcal{C}_f$ as (\ref{param}), and $\omega=A(x_1,x_2)dx_1+B(x_1,x_2)dx_2\in\Omega_f$, then we denote by $$\varphi^*(\omega):=A(\varphi(T))\cdot\varphi'_1(T)+B(\varphi(T))\cdot\varphi'_2(T)$$ the pullback of $\omega$ by $\varphi$, where $\varphi'_i(T)$ denotes the derivative of $\varphi_i(T)$ with respect to $T$, for $i=1,2$. 

The \emph{value of a $1$-form} $\omega\in\Omega_f$ (with respect to $\mathcal{C}_f$) is $\nu_f(\omega):=ord_T\left ( \varphi^*(\omega)\right )+1$ and we define
$$\Lambda_f:=\nu_f(\Omega_f):=\{\nu_f(\omega);\ \omega\in\Omega_f\ \mbox{such that}\ \varphi^*(\omega)\neq 0\}$$
the \emph{value set of $1$-form} with respect to $\mathcal{C}_f$.

The set $\Lambda_f$ is an analytic invariant of the plane branch $\mathcal{C}_f$ and can be computed using the algorithm presented in \cite{basestandard}. In \cite{classification} (and \cite{handbook}), $\Lambda_f$ plays a central role in the solution to the analytic classification problem for plane branches. 

Carbonne \cite{carbonne} defines an intermediate equivalence between the topological and the analytical equivalence, called {\it equidifferentiable equivalence}: two plane branches $f$ and $h$ are equidifferentiable if $\Gamma_f=\Gamma_h$ and $\Lambda_f=\Lambda_h$. 

Note that $\nu_f(h\omega)=\nu_f(h)+\nu_f(\omega)$, and if $0\neq h\in\langle X_1,X_2\rangle\subset\mathcal{O}_2$ then $\nu_f(h)=\nu_f(dh)$. Thus, we obtain $\Gamma_f+\Lambda_f\subseteq \Lambda_f$ and $\Gamma_f\setminus\{0\}\subseteq\Lambda_f$. In \cite{torsion}, Zariski showed that $\Lambda_f=\Gamma_f\setminus\{0\}$ if and only if $\mathcal{C}_f$ is analytically equivalent to $\mathcal{C}_h$, where $h=X_2^{v_0}-X_1^{v_1}$ and $\Gamma_f=\langle v_0,v_1\rangle$. 

If $\Lambda_f\neq\Gamma_f\setminus\{0\}$ and $\lambda_1:=\min\Lambda_f\setminus\Gamma_f$, then $\lambda_0:=\lambda_1-v_0$ is called {\it Zariski invariant} of $\mathcal{C}_f$. According to \cite{classification}, there exists a plane branch analytically equivalent to $\mathcal{C}_f$ that admits a parameterization of the form \begin{equation}\label{forma-normal}\left ( T^{v_0}, T^{v_1}+a_{\lambda_0}T^{\lambda_0}+\sum_{i+v_0\not\in\Lambda_f} a_iT^{i}\right )\ \ \mbox{with}\ \ a_{\lambda_0}\neq 0.\end{equation}
Note that, for the above parameterization we get
$\nu_f\left (\frac{1}{v_0}x_2dx_1-\frac{1}{v_1}x_1dx_2\right )=\lambda_0+v_0=\lambda_1$.

\begin{remark}\label{obs-zariski}
Notice that $v_1<\lambda_0$, and $\lambda_0\leq\beta_2$ if $g\geq 2$, where $\beta_2$ is as (\ref{expoentes-caract}). Thus, $\lambda_1=\min\Lambda_f\setminus\Gamma_f$ satisfies $\lambda_1>v_1+v_0$. Moreover, 
Zariski, in \cite{torsion}, showed that there is no $\lambda\in\Lambda_f\setminus\Gamma_f$ such that $\lambda_1<\lambda<\lambda_1+v_0$.
\end{remark}

The set $\Lambda_f$ encodes the generators of the value semigroup; that is, we can recover $\Gamma_f$ from $\Lambda_f$ (see \cite{osnar}). In addition, $\Lambda_f$ is related to the Tjurina number $\tau_f$ of $f$. More precisely, we have (cf. \cite[Proposition 1.2.19]{handbook})
\begin{equation}\label{tjurina-lambda}
\tau_f=\mu_f-\sharp\Lambda_f\setminus\Gamma_f .
\end{equation}

\begin{example}\label{6-9-v2}
For any plane branch with value semigroup $\Gamma=\langle 6,9,v_2\rangle$, we have $\dim\frac{\widetilde{H''}}{H''}\geq 5$.

In fact, since $\dim\frac{\widetilde{H''}}{H''}\geq \mu-\tau=\sharp\Lambda\setminus\Gamma$, it suffices to show that for any such plane branch, $\sharp\Lambda\setminus\Gamma\geq 5$.

According to Theorem 1.2.7 and Theorem 1.2.25 in \cite{handbook}, any plane branch $\mathcal{C}_f$ with value semigroup $\Gamma=\langle 6,9,v_2\rangle$ admits (up to analytical equivalence) a parameterization of the form
$$\varphi (T)=\left (T^6,T^9+T^{v_2-9}+\sum_{j=-\left [ \frac{v_2}{3}\right ]-3}^{-4}a_{2v_2+3j}T^{2v_2+3j} \right ).$$

Considering $h=X_2^2-X_1^3$, we obtain:
$$h(\varphi(T))=2T^{v_2}+\sum_{j=-\left [ \frac{v_2}{3}\right ]-3}^{-10}2a_{2v_2+3j}T^{2v_2+3(j+3)}+(1+2a_{2v_2-27})T^{2v_2-18}+\mbox{h.o.t.}\footnote{Here h.o.t. means higher order terms.}$$ where $a_{2v_2+3j}=0$ for $-\left [ \frac{v_2}{3}\right ]-3\leq j\leq -10$ if $v_2<21$.

Taking $\omega_0=3x_2dx_1-2x_1dx_2$ we get:
$$\nu_f(\omega_0)=v_2-3,\ \ \ \nu_f(x_1\omega)=v_2+3,\ \ \ \nu_f(h\omega_0)=2v_2-3\ \ \ \mbox{and}\ \ \ \nu_f(x_1h\omega_0)=2v_2+3,$$
that is, $\{v_2-3, v_2+3, 2v_2-3, 2v_2+3\}\subseteq\Lambda\setminus\Gamma$.

Moreover, consider
{\small
$$T\cdot\varphi^*\left (\frac{6x_1dh-v_2hdx_1}{12}\right )=\sum_{j=-\left [ \frac{v_2}{3}\right ]-3}^{-10}(v_2+3(j+3))a_{2v_2+3j}T^{2v_2+3(j+5)}+\frac{(v_2-18)}{2}(2a_{2v_2-27}+1)T^{2v_2-12}+\mbox{h.o.t.}$$}
{\small $$T\cdot\varphi^*\left (\frac{9x_2dh-v_2hdx_2}{18}\right )=\sum_{j=-\left [ \frac{v_2}{3}\right ]-3}^{-10}(v_2+3(j+3))a_{2v_2+3j}T^{2v_2+3(j+6)}+\frac{(v_2-18)}{18}(18a_{2v_2-27}+9-2v_2)T^{2v_2-9}+\mbox{h.o.t.}.$$
}

If there exists $a_{2v_2+3j}\neq 0$ for some $-\left [ \frac{v_2}{3}\right ]-3\leq j\leq -10$, then 
$$\{v_2-3, v_2+3, 2v_2-3, 2v_2+3,\nu_f(6x_1dh-v_2hdx_1), \nu_f(9x_2dh-v_2hdx_2)\}\subseteq\Lambda\setminus\Gamma,$$ so $\Lambda\setminus\Gamma\geq 6$.

If $a_{2v_2+3j}= 0$ for all $-\left [ \frac{v_2}{3}\right ]-3\leq j\leq -10$, then 
$$\varphi^*\left (\frac{6x_1dh-v_2hdx_1}{12}\right )=\frac{(v_2-18)}{2}(2a_{2v_2-27}+1)T^{2v_2-12}+\mbox{h.o.t.}\ \ \mbox{and}$$
$$\varphi^*\left (\frac{9x_2dh-v_2hdx_2}{18}\right )=\frac{(v_2-18)}{18}(18a_{2v_2-27}+9-2v_2)T^{2v_2-9}+\mbox{h.o.t.}.$$ Thus 
$$\{v_2-3, v_2+3, 2v_2-3, 2v_2+3,2v_2-9\}\subseteq\Lambda\setminus\Gamma\ \ \mbox{or}\ \ \{v_2-3, v_2+3, 2v_2-3, 2v_2+3,2v_2-12\}\subseteq\Lambda\setminus\Gamma,$$ and $\Lambda\setminus\Gamma\geq 5$.

In \cite{6-9-19}, all possible sets $\Lambda\setminus\Gamma$ for $\Gamma=\langle 6,9,19\rangle$ are presented, and it is shown that $\sharp\Lambda\setminus\Gamma=5$ can indeed be achieved for some plane branches with value semigroup $\Gamma$.

Hence, for any plane curve with value semigroup $\langle 6,9,v_2\rangle$ we have $\dim\frac{\widetilde{H''}}{H''}\geq \sharp\Lambda\setminus\Gamma\geq 5$.
\end{example}

\section{The sets $\Lambda_f$ and $\mathcal{R}_f$ for $f\in\mathcal{O}_2$ with $\Gamma_f=\langle v_0,v_1\rangle$}

For any irreducible $f\in\mathcal{O}_2$, by Theorem \ref{TheoremX}, (\ref{mu-tau}) and (\ref{tjurina-lambda}), we have that for each $\lambda\in\Lambda_f\setminus\Gamma_f$, there exists $\sigma_1+1<\sigma_{\lambda}\in\mathcal{S}_f$ such that $-\sigma_{\lambda}\in\mathcal{R}_f$. In this section, we show how to determine $\sigma_{\lambda}\in\mathcal{S}_f$, and consequently $-\sigma_{\lambda}\in\mathcal{R}_f$, for any $\lambda\in\Lambda_f\setminus\Gamma_f$ when $f\in\mathcal{O}_2$ is semiquasihomogeneous, that is, when $\Gamma_f=\langle v_0,v_1\rangle$.

According to Zariski (cf. \cite[Chapter VI]{zariski}), up to an analytic change of coordinates, any plane branch with $\Gamma=\langle v_0,v_1\rangle$ can be defined by 
\begin{equation}\label{eqf}
f=X_2^{v_0}-X_1^{v_1}+\sum_{(i,j)\in\Theta}a_{ij}X_1^iX_2^j\end{equation}
where $\Theta=\{(i,j)\in\mathbb{N}^2;\ 1\leq i<v_1-1,\ 1\leq j<v_0-1\ \mbox{and}\ iv_0+jv_1>v_0v_1\}$.

Since any $f\in\mathcal{O}_2$ of the form (\ref{eqf}) is Newton nondegenerate with respect to its Newton polygon, we obtain an alternative way to compute the $\mathcal{V}$-filtration on $H''/H'$, which we describe in the sequel.

Consider the (local) monomial order $\preceq$ on $\mathcal{O}_2$ defined by
$$X_1^{i_1}X_2^{j_1}\prec X_1^{i_2}X_2^{j_2}\ \Leftrightarrow\ \left \{ \begin{array}{l}i_1v_0+j_1v_1<i_2v_0+j_2v_1\ \ \mbox{or}\\
i_1v_0+j_1v_1=i_2v_0+j_2v_1\ \mbox{and}\ X_1^{i_1}X_2^{j_1}<_{lex} X_1^{i_2}X_2^{j_2}\end{array}\right .$$
where $\leq_{lex}$ denotes the lexicographic order with $X_1\leq_{lex} X_2$.

Since $\{f_{X_1}, f_{X_2}\}$ is a Standard basis for the Jacobian ideal $J(f)=(f_{X_1},f_{X_2})$ with respect to $\preceq$, it follows that
\begin{equation}\label{E}E=\{\overline{X}_1^i\overline{X}_2^j;\ 0\leq i< v_1-1\ \mbox{and}\ 0\leq j<v_0-1\}\end{equation}
is a $\mathbb{C}$-basis of $\mathcal{O}_2/J(f)$. 
Given $\overline{0}\neq\overline{h}=\sum_{i,j}c_{ij}\overline{X}_1^i\overline{X}_2^j\in\frac{\mathcal{O}_2}{J(f)}$, the {\it Newton order} of $\overline{h}$ (with respect to $f$) is defined as
\[N_f(\overline{h})=\inf\left \{ \frac{i}{v_1}+\frac{j}{v_0};\ c_{ij}\neq 0\right \}.\]

According to Varchenko and Khovanskii (cf. \cite{VK}), the filtration induced by the Newton order on $\mathcal{O}_2/J(f)$ and the $\mathcal{V}$-filtration are related by 
\begin{equation}\label{equiv}
\overline{[hdX_1\wedge dX_2]}\in \text{Gr}_{\mathcal{V}}^{\alpha}(H''/H')\ \ \mbox{if and only if}\ \ \alpha=N_f(\overline{X_1X_2h})-1.\end{equation}

Notice that if $\overline{0}\neq\overline{h}=\sum_{i,j}c_{ij}\overline{X}_1^i\overline{X}_2^j\in\mathcal{O}_2/J(f)$ then
\begin{equation}\label{newton-filtration}
N_f(\overline{X_1X_2h})-1 = \inf\left \{ \frac{i+1}{v_0}+\frac{j+1}{v_1}-1;\ c_{ij}\neq 0\right \} =\frac{\inf \{ iv_0+jv_1;\ c_{ij}\neq 0 \}-(\mu_f-1)}{v_0v_1}, 
\end{equation}
where $\mu_f=(v_0-1)(v_1-1)$ is the Milnor number of $f$, which coincides with the conductor $c_f$ of $\Gamma_f=\langle v_0,v_1\rangle$.
Moreover, according to Saito (cf. \cite{Saito-93}), we have that
\begin{equation}\label{S}
\mathcal{S}_f=\left \{N_f\left (\overline{X}_1^{i+1}\overline{X}_2^{j+1}\right );\ \overline{X}_1^i\overline{X}_2^j\in E \right \}=\left \{ \frac{i+1}{v_0}+\frac{j+1}{v_1}-1;\ 0\leq i<v_1-1\ \mbox{and}\ 0\leq j<v_0-1\right \}.\end{equation}
In particular, $\sigma_1=\min \mathcal{S}_f=-\frac{\mu_f-1}{v_0v_1}$.

It follows from (\ref{eqf}) that 
\begin{equation}\label{f-equat}
f=\frac{1}{v_1}X_1f_{X_1}+\frac{1}{v_0}X_2f_{X_2}-\sum_{(i,j)\in\Theta}\left( \frac{i}{v_0}+\frac{j}{v_1} -1\right )a_{ij}X_1^iX_2^j\end{equation}
and any $\overline{0}\neq\overline{h}\in\frac{(f)+J(f)}{J(f)}$ can be expressed as $\overline{h}=\sum_{(i,j)\in\Theta}c_{ij}\overline{X}_1^i\overline{X}_2^j$. Notice that $\nu_f(X_1^iX_2^j)=iv_0+jv_1$ and $\nu_f(X_1^{i_1}X_2^{j_1})\neq\nu_f(X_1^{i_2}X_2^{j_2})$ for every distinct elements $(i_1,j_1)$ and $(i_2,j_2)$ in $\Theta$, thus by (\ref{newton-filtration}):
$$N_f(\overline{X_1X_2h})-1=\frac{\nu_f(\sum_{(i,j)\in\Theta}c_{ij}X_1^iX_2^j)-(\mu_f-1)}{v_0v_1}.$$

It follows from (\ref{E}) that  $$\nu_f\left (\frac{J(f)}{(f)}\right )=\mu_f-1+\Gamma\setminus\{0\}$$
and, according to \cite[Proposition 3.31]{Pol} and \cite[Theorem 5]{tjurina}, we have
$$\nu_f\left (\frac{(f)+J(f)}{(f)}\right )=\mu_f-1+\nu_f(\Omega_f)=\mu_f-1+\Lambda_f.$$
Moreover, if $Af_{X_1}+Bf_{X_2}\in \frac{(f)+J(f)}{(f)}$, then setting $Bdx_1-Adx_2\in\Omega_f$ we obtain
$$\nu_f(Af_{X_1}+Bf_{X_2})=\mu_f-1+\nu_f(Bdx_1-Adx_2).$$

In this way, we have the following:

\begin{theorem}\label{TheoremXX} Let $C_f$ be a plane branch with value semigroup $\Gamma_f=\langle v_0,v_1\rangle$ and value set of $1$-form $\Lambda_f$. For every $\lambda\in\Lambda_f\setminus\Gamma_f$, we have that $-\frac{\lambda}{v_0v_1}\in \mathcal{R}_f$.
\end{theorem}
\Dem Up to analytic change of coordinates, we may assume that $f$ is given by (\ref{eqf}). Since,  $\nu_f\left (\frac{(f)+J(f)}{(f)}\right )=\mu_f-1+\nu_f(\Omega_f)=\mu_f-1+\Lambda_f$ and $\nu_f\left (\frac{J(f)}{(f)}\right )=\mu_f-1+\Gamma_f\setminus\{0\}$, it follows that for any $\lambda\in\Lambda_f\setminus\Gamma_f$, there exists an element $h\in (f)+J(f)\subset\mathcal{O}_2$ such that 
$\overline{0}\neq\overline{h}=\sum_{(i,j)\in\Theta}c_{ij}\overline{X}_1^i\overline{X}_2^j\in\frac{(f)+J(f)}{J(f)}$ and
$$\lambda+\mu_f-1=\nu_f(h)=\nu_f\left ( \sum_{(i,j)\in\Theta}c_{ij}X_1^iX_2^j\right )=\inf\{iv_0+jv_1;\ c_{ij}\neq 0\}.$$

By (\ref{equiv}) and (\ref{newton-filtration}), it follows that $N_f(\overline{X_1X_2h})-1=\frac{v_f(h)-(\mu_f-1)}{v_0v_1}=\frac{\lambda}{v_0v_1}$, and consequently,
$[\overline{hdX_1\wedge dX_2}]\in \text{Gr}_{\mathcal{V}}^{\frac{\lambda}{v_0v_1}}(H''/H')$. Then, by Theorem \ref{TheoremX}, we conclude that $-\frac{\lambda}{v_0v_1}\in\mathcal{R}_f$.
\cqd

Notice that the previous theorem allows us to relate spectral numbers with roots of $\tilde{b}_f$ by means of the analytical invariant $\Lambda_f\setminus\Gamma_f$ of $\mathcal{C}_f$. In particular, Theorem \ref{TheoremXX} provides a proof of Conjecture 1.3 in \cite{david}.

\begin{lemma}\label{espectral-lacunas}
For any plane branch $\mathcal{C}_f$ with $\Gamma_f=\langle v_0,v_1\rangle$ we have $$\mathcal{S}_f=\left \{ -\frac{\delta}{v_0v_1}, \frac{\delta}{v_0v_1};\ 0<\delta\in\mathbb{N}\setminus\Gamma_f\right \}.$$
\end{lemma}
\Dem
By Remark \ref{escrita}, we have $\{\delta;\ 0<\delta\not\in\Gamma_f\}=\{\alpha v_1-\beta v_0;\ 1\leq \alpha<v_0\ \mbox{and}\ 1\leq \beta<\frac{\alpha v_0}{v_1}\}$.

On the other hand, by (\ref{S}), given $\sigma\in \mathcal{S}_f$ we have 
$$\sigma=\frac{(i+1)v_0+(j+1)v_1-v_0v_1}{v_0v_1}\ \ \mbox{with}\ \ 0\leq i<v_1-1\ \ \mbox{and}\ \ 0\leq j<v_0-1.$$

If $i<\frac{(v_0-1-j)v_1}{v_0}-1$, then
$$(i+1)v_0+(j+1)v_1-v_0v_1=-((v_0-1-j)v_1-(i+1)v_0)$$
with $1\leq v_0-1-j<v_0$ and $1\leq i+1<\frac{(v_0-1-j)v_1}{v_0}.$
Hence, $$\{-((v_0-1-j)v_1-(i+1)v_0);\ 1\leq v_0-1-j<v_0\ \ \mbox{and}\ \ 1\leq i+1<\frac{(v_0-1-j)v_1}{v_0}\}=\{-\delta;\ 0<\delta\not\in\Gamma_f\}.$$

If $i>\frac{(v_0-1-j)v_1}{v_0}-1$, then
$$(i+1)v_0+(j+1)v_1-v_0v_1=(j+1)v_1-(v_1-1-i)v_0\ \ \mbox{with}\ \ 1\leq j+1<v_0\ \ \mbox{and}\ \ 1\leq v_1-1-i<\frac{(j+1)v_1}{v_0}.$$
Thus, $\{(j+1)v_1-(v_1-1-i)v_0;\ 1\leq j+1<v_0\ \ \mbox{and}\ \ 1\leq v_1-1-i<\frac{(j+1)v_1}{v_0}\}=\{\delta;\ 0<\delta\not\in\Gamma_f\}$.

Hence, $\mathcal{S}_f=\left \{ -\frac{\delta}{v_0v_1}, \frac{\delta}{v_0v_1};\ 0<\delta\not\in\Gamma_f\right \}.$
\cqd

As a consequence of the previous lemma and Theorem \ref{TheoremXX} we have the following result.

\begin{corollary}\label{condutor}
If $\Lambda_f\setminus\Gamma_f\neq \emptyset$ and $\lambda_c:=\min\{\delta\in\Lambda_f\setminus\Gamma_f;\ \gamma\in\Lambda_f\ \mbox{for every}\ \gamma\geq\delta\}$, then $-\sigma\in\mathcal{R}_f$ for every $\frac{\lambda_c}{v_0v_1}\leq\sigma\in\mathcal{S}_f$.    
\end{corollary}
\Dem According to Lemma \ref{espectral-lacunas}, given $\sigma\in\mathcal{S}_f$, we have $\sigma=-\frac{\gamma}{v_0v_1}$ or $\sigma=\frac{\gamma}{v_0v_1}$ with $\gamma\not\in\Gamma_f$. If $\frac{\lambda_c}{v_0v_1}\leq\sigma=\frac{\gamma}{v_0v_1}$, i.e., $\gamma\geq\lambda_c$, then $\gamma\in\Lambda_f\setminus\Gamma_f$. So, by Theorem \ref{TheoremXX}, we obtain $-\sigma=-\frac{\gamma}{v_0v_1}\in\mathcal{R}_f$.
\cqd

Let us assume that $f\in\mathcal{O}_2$ is given as (\ref{f-equat}). Thus, we have
\begin{equation}\label{f1}
f_1:=f-\frac{1}{v_1}X_1f_{X_1}-\frac{1}{v_0}X_2f_{X_2}=-\sum_{(i,j)\in\Theta}\left( \frac{i}{v_0}+\frac{j}{v_1} -1\right )a_{ij}X_1^iX_2^j\in (f)+J(f)\end{equation}
and (cf. \cite{KSaito}) $f_1\in J(f)$ if and only if the branch $\mathcal{C}_f$ is analytically equivalent to the branch defined by $X_2^{v_0}-X_1^{v_1}$, that is, $\Lambda_f=\Gamma_f\setminus\{0\}$, or equivalently $\tau_f=\mu_f$.

Now, assuming $\Lambda_f\neq\Gamma_f\setminus\{0\}$, a parameterization of $\mathcal{C}_f$ is given by
$$\left ( T^{v_0},T^{v_1}-\frac{a_{i_0j_0}}{v_0}T^{\lambda_0}+h.o.t.\right ),$$
where $i_0v_0+j_0v_1=\inf\{iv_0+jv_1;\ a_{ij}\neq 0\}$, and $\lambda_0=i_0v_0+j_0v_1-(v_0-1)v_1$ is the Zariski invariant of $\mathcal{C}_f$, that is,
$$\lambda_1=\lambda_0+v_0=\inf\{iv_0+jv_1;\ a_{ij}\neq 0\}-(\mu_f-1)=\min\Lambda_f\setminus\Gamma_f.$$ In this way, we obtain $\overline{0}\neq\overline{f}_1\in\frac{(f)+J(f)}{J(f)}$, and
$$\inf\{iv_0+jv_1;\ a_{ij}\neq 0\}=\nu_f(f_1)=\mu_f-1+\nu_f\left (\frac{1}{v_0}x_2dx_1-\frac{1}{v_1}x_1dx_2\right )=\mu_f-1+\lambda_1.$$
Consequently, by (\ref{newton-filtration}), we have
$$N_f(\overline{X_1X_2f_1})-1=\frac{\inf\{iv_0+jv_1;\ a_{ij}\neq 0\}-(\mu_f-1)}{v_0v_1}=\frac{\lambda_1}{v_0v_1}$$
and by (\ref{equiv}) we obtain $[\overline{f_1dX_1\wedge dX_2}]\in \text{Gr}_{\mathcal{V}}^{\frac{\lambda_1}{v_0v_1}}\left ( \frac{H''}{H'}\right )$ and $[f_1dX_1\wedge dX_2]\in \mathcal{V}^{\frac{\lambda_1}{v_0v_1}}H''\subseteq  \mathcal{V}^{\frac{\lambda_1}{v_0v_1}}\widetilde{H''}$.

In this way, we obtain the following result.

\begin{proposition}\label{limitante-inf} Given a plane branch $\mathcal{C}_f$ with value semigroup $\Gamma_f=\langle v_0,v_1\rangle$ and value set of $1$-form $\Lambda_f$, if $\lambda_1=\min(\Lambda_f\setminus\Gamma_f)$, then $-(\sigma+1)\in\mathcal{R}_f$ for every $\sigma\in \mathcal{S}_f$ such that $\sigma<\frac{\lambda_1}{v_0v_1}$.
\end{proposition}
\Dem Let us consider $f_1\in (f)+J(f)$ as (\ref{f1}). Since $f-f_1\in J(f)$ and $[f_1dX_1\wedge dX_2]\in \mathcal{V}^{\frac{\lambda_1}{v_0v_1}}\widetilde{H''}$, then by \cite[Proposition 3.1]{andrea} we have $\widetilde{H'}\subseteq H'+\mathcal{V}^{\frac{\lambda_1}{v_0v_1}}\widetilde{H''}$, and by \cite[Proposition 3.2]{andrea}, it follows that for every spectral number $\sigma\in\mathcal{S}_f$ with $\sigma<\frac{\lambda_1}{v_0v_1}$, the number $-(\sigma+1)$ is a root of $\widetilde{b}_f$. 
\cqd

Notice that, by setting $\lambda_1=\infty$ when $\Lambda_f\setminus\Gamma_f=\emptyset$, the previous proposition recovers the classical relation $\mathcal{R}_f=-(\mathcal{S}_f+1)$ if $f$ is analytically equivalent to a quasihomogeneous element of $\mathcal{O}_2$. 

Let us consider the following partition of the spectral numbers $\mathcal{S}_f$ associated to a plane branch $\mathcal{C}_f$ with semigroup $\Gamma_f=\langle v_0,v_1\rangle$:
$$\mathcal{S}_f=\mathcal{S}_{-1}\ \cup\ \widetilde{\mathcal{S}}_f\ \cup\ \mathcal{S}_0$$
where
$$\mathcal{S}_{-1}=\left \{\sigma\in\mathcal{S}_f;\  \sigma<\frac{\lambda_1}{v_0v_1}\right \},
\ \  \widetilde{\mathcal{S}}_f=\left \{\sigma=\frac{\gamma}{v_0v_1}\in\mathcal{S}_f;\ \lambda_1<\gamma\not\in\Lambda_f\right \},\ \ \
\mathcal{S}_0=\left \{\sigma=\frac{\lambda}{v_0v_1}\in\mathcal{S}_f;\  \lambda\in\Lambda_f\right \}.$$
According to Theorem \ref{TheoremXX} and Proposition \ref{limitante-inf}, we have 
$$\{-(\sigma+1);\ \sigma\in\mathcal{S}_{-1}\}\ \cup\ \{-\sigma;\ \sigma\in\mathcal{S}_0\}\ \subset\ \mathcal{R}_f,$$
that is, by using the value set $\Lambda_f$, or equivalently, fixing an equidifferentiable class of plane branches, we can relate a subset of the spectral numbers to roots of the Bernstein polynomial of $f\in\mathcal{O}_2$. Note that for $\sigma \in\widetilde{\mathcal{S}}_f$, we may have either $-(\sigma+1)\in\mathcal{R}_f$ or $-\sigma\in\mathcal{R}_f$.

In particular, if $\mathcal{G}_f=\{l>\lambda_1;\ l\not\in\Lambda_f\}$ then it follows, from (\ref{mu-tau}) and (\ref{tjurina-lambda}), that
\begin{equation}\label{desig}
\sharp\Lambda_f\setminus\Gamma_f\leq
\dim_{\mathbb{C}}\frac{\widetilde{H''}}{H''}=\sum_{\sigma\in\mathcal{S}_f}\sigma-\sum_{\epsilon\in\mathcal{E}_f}\epsilon\leq\sharp\Lambda_f\setminus\Gamma_f+\sharp\widetilde{\mathcal{S}}_f=\sharp\Lambda_f\setminus\Gamma_f+\sharp\mathcal{G}_f.\end{equation}

For any plane branch $\mathcal{C}_f$, not necessarily with semigroup generated by two elements, we have that $\mathcal{G}_f=\{l>\lambda_1;\ l\not\in\Lambda_f\}=\emptyset$ if and only if every integer greater than or equal to $\lambda_1$ belongs to $\Lambda_f$. Equivalently, this occurs if $\min\Lambda_f\setminus\Gamma_f=\lambda_1=\lambda_c=\min\{\lambda\in\Lambda_f;\ \lambda+\mathbb{N}\subset\Lambda_f\}$ and consequently, $\Lambda_f\setminus\Gamma_f$ is completely determined by $\Gamma_f$ and $\lambda_1$.

The following result characterizes all plane branches $\mathcal{C}_f$ such that $\mathcal{G}_f=\emptyset$.

\begin{proposition}\label{G=0}
We have $\mathcal{G}=\{l>\lambda_1;\ l\not\in\Lambda\}=\emptyset$ for a plane branch $\mathcal{C}$ if and only if: 
\begin{enumerate}
\item $\Gamma=\langle v_0,v_1\rangle$ and $\lambda_1=\infty$ or;
\item $\Gamma=\langle v_0,v_1\rangle$ and $\lambda_1=(v_0-1)v_1-sv_0$ for some $1\leq s\leq \left [\frac{v_1}{v_0} \right ]+1$ or;
\item $\Gamma=\langle 4,6,v_2\rangle$.
\end{enumerate}
\end{proposition}
\Dem It is clear that for any plane curve analytically equivalent to the monomial curve defined by $X_2^{v_0}-X_1^{v_1}$, we have $\lambda_1=\infty$ and $\Lambda\setminus\Gamma=\emptyset$. Hence, $\mathcal{G}=\emptyset$.	

Thus, it is sufficient to consider the case $\lambda_1<\infty$, that is, $\Lambda\setminus\Gamma\neq\emptyset$. 
	
Notice that if $\Gamma=\langle v_0,v_1,\ldots ,v_g\rangle$ with $g\geq 2$ then, by Remark \ref{obs-zariski} and (\ref{exp-ger}) we have
$$\lambda_1\leq \beta_2+v_0=\left \{\begin{array}{ll}v_2+v_1-\left (\frac{v_1}{e_1}-1\right )v_0<v_3-v_0 & \mbox{if}\ g\geq 3 \\
v_2-(n_1-1)v_1+v_0 < v_2-v_0 & \mbox{if}\ g=2\ \mbox{and}\ n_1\geq 3 \\
v_2-v_1+v_0 < v_2-v_0 & \mbox{if}\ g=2, n_1=2\ \mbox{and}\ v_1>2v_0.\end{array}\right .$$

By Lemma 2.6 in \cite{osnar}, we have $v_i-v_0\not\in\Lambda$ for any plane curve with semigroup $\Gamma$ and $1\leq i\leq g$. In this way, in the above cases, we have $\mathcal{G}\neq\emptyset$. Therefore, in order to have $\mathcal{G}=\emptyset$, we must have: 
		
\begin{enumerate}
\item[\mbox{\bf Case i:}] $g=2$, $n_1=2$ and $v_1<2v_0$. 

In this case, by (\ref{exp-ger}), we have $v_2=v_1+\beta_2$ and
$2e_1=n_1e_1=v_0<v_1<2v_0=2n_1e_1=4e_1$, which implies $\frac{v_1}{e_1}=3$, and consequently $\Gamma=\langle 2e_1, 3e_1, v_2\rangle$.
		
We must have $\lambda_1=\beta_2+v_0$, since otherwise, by Remark \ref{escrita}, we would obtain $\lambda_1=s_1v_1+s_0v_0$ with $0\leq s_1<n_1=2$ and $s_0<0$, implying $\lambda_1<v_1$, which is absurd. Since $2v_1=3v_0$ and $v_2=v_1+\beta_2$, it follows that $\lambda_1=\beta_2+v_0=v_2-v_1+v_0=v_2+v_1-2v_0.$ Let $0<k\in\mathbb{N}$ such that $1<k-\frac{\beta_2}{v_0}<2$. Then
$$\lambda_1=v_2+v_1-2v_0<2v_2-kv_0<v_2+v_1-v_0=\lambda_1+v_0.$$
If $n_2>2$, then $2v_2-kv_0\not\in\Gamma$ so, by Remark \ref{obs-zariski}, we have $2v_2-kv_0\in \mathcal{G}$. 
		
Thus, to ensure $\mathcal{G}=\emptyset$, we must have $e_1=n_2=2$, which gives $\Gamma_f=\langle 4, 6, v_2\rangle$, and consequently $\lambda_1=\beta_2+v_0=v_2-v_1+v_0=v_2-2$. Since the conductor of the semigroup $\Gamma$ is given by $\mu=(n_2-1)v_2+(n_1-1)v_1-v_0+1=v_2+3$, the even numbers $v_2-1$ and $v_2+1$ belong to $\langle 4,6\rangle\subset\Gamma$. Thus, $\{\lambda_1=v_2-2,v_2-1,v_2,v_2+1,\lambda_1+v_0=v_2+2,v_2+3,\ldots\}\subset\Lambda$. In particular, $\Lambda\setminus\Gamma=\{v_2-2, v_2+2\}$ and $\mathcal{G}=\emptyset$.
		
\item[\mbox{\bf Case ii:}] $g=1$ and $\lambda_1<\infty$.
		
Since $\Gamma=\langle v_0,v_1\rangle$ and $v_1+v_0<\lambda_1\not\in\Gamma$ we may write $\lambda_1=rv_1-sv_0$ with $s>0$ and $2\leq r\leq v_0-1$. 
		
Notice that $\lambda_1+\Gamma\subset\Lambda$, so $(v_0-1)v_1-sv_0\in\Lambda$. 		
If $r<v_0-1$, take $k=\left [ \frac{(v_0-1-r)v_1}{v_0}\right ]+s$, then
$$\lambda_1=rv_1-sv_0<(v_0-1)v_1-kv_0<rv_1-(s-1)v_0=\lambda_1+v_0.$$
As $k>0$, we have $(v_0-1)v_1-kv_0\not\in\Gamma$, so by Remark \ref{obs-zariski}, it follows that $(v_0-1)v_1-kv_0\in\mathcal{G}$. Therefore, in order to have $\mathcal{G}=\emptyset$, we must have $\lambda_1=(v_0-1)v_1-sv_0$ for some $s\geq 1$, and $\{\lambda_1+nv_0;\ 0\leq n<s-1\}=\{(v_0-1)v_1-jv_0;\ 1\leq j\leq s\}\subseteq\Lambda\setminus\Gamma$.
		
If $s>\frac{v_1}{v_0}+1$, take $k\in\mathbb{N}$ such that $0<s-\frac{v_1}{v_0}-k<1$. Then
$$\lambda_1=(v_0-1)v_1-sv_0<(v_0-2)v_1-kv_0<(v_0-1)v_1-sv_0+v_0=\lambda_1+v_0.$$
Since $(v_0-2)v_1-kn\not\in\Gamma$, again by Remark \ref{obs-zariski}, we have $(v_0-1)v_1-kv_0\in\mathcal{G}$.
		
Therefore, we must have $\lambda_1=(v_0-1)v_1-sv_0$ for some $1\leq s\leq \left [ \frac{v_1}{v_0}\right ]+1$.

Now, consider $\ell=r_1v_1+r_0v_0>\lambda_1$ with $0\leq r_1<v_0$. If $r_1=v_0-1$, then $r_1>-s$ and $\ell\in\Lambda$. If $r_1<v_0-1$, then $r_0v_0>(v_0-1-r_1)v_1-sv_0>-v_0$, i.e., $r_0\geq 0$ and $\ell\in\Gamma\subset\Lambda$. Hence,  $\mathcal{G}=\emptyset$.
\end{enumerate}
\cqd

As recalled in Section \ref{sec:plana}, for any plane branch in the topological class determined by the value semigroup $\langle 4,6,v_2\rangle$, we have $\dim\frac{\widetilde{H''}}{H''}=\mu-\tau=\sharp\Lambda\setminus\Gamma=2$. This property, along with the explicit description of $\mathcal{R}$, was obtained by Bartolo, Cassous-Nogu\`es, Luengo and Melle-Hern\'andez (see Theorem 4.2 in \cite{cassous-2}).

In the other cases described in Proposition \ref{G=0}, the topological class is determined by $\Gamma=\langle v_0,v_1\rangle$, so by (\ref{desig}), we again have $\dim\frac{\widetilde{H''}}{H''}=\mu-\tau=\sharp\Lambda\setminus\Gamma$. Note that among these cases,
we find all plane branches with $2\leq v_0\leq 3$, all plane branches such that $\tau=\mu-1$ (see \cite{bayer} or Corollary 1.4.4 in \cite{handbook}), and all plane branches with $0$-dimensional moduli space (see \cite{gaffney} or Proposition 1.4.5 in \cite{handbook}).
	
In the sequel, we present analytical invariants for plane branches satisfying the above proposition.
	
\begin{center}
\begin{tabular}{|c|c|c|c|}
\hline
{\bf $\Gamma,\ \mu,\ \lambda_1$} & {\bf $\Lambda\setminus\Gamma$} &  {\bf $\tau=\mu-\sharp\Lambda\setminus\Gamma$} & {$\mathcal{R}$}\\
\hline

\footnotesize{$\begin{array}{c}\Gamma=\langle v_0, v_1\rangle \\ \mu=(v_0-1)(v_1-1)\\ \lambda_1=\infty\end{array}$} & $\emptyset$ & \footnotesize{$(v_0-1)(v_1-1)$} & \footnotesize{$\begin{array}{c}
-\left ( \frac{l}{v_0v_1}+1\right ),\ 	-\left ( \frac{-l}{v_0v_1}+1\right )\\ 0<l\not\in\Gamma\end{array}$}\\ \hline

\footnotesize{$\begin{array}{c}\Gamma=\langle v_0, v_1\rangle \\ \mu=(v_0-1)(v_1-1)\\ \lambda_1=(v_0-1)v_1-sv_0 \\ 1\leq s\leq\left [ \frac{v_1}{v_0}\right ]+1\end{array}$} & \footnotesize{$\begin{array}{c} \lambda_1+jv_0 \\ 0\leq j\leq s-1\end{array}$} &  \footnotesize{$(v_0-1)(v_1-1)-s$} & \footnotesize{$\begin{array}{c}
-\left ( \frac{l_0}{v_0v_1}+1\right ),\ -\left ( \frac{l_1}{v_0v_1}+1\right ),\\	-\frac{(v_0-1)v_1-jv_0}{v_0v_1}\\ 0<l_0<(v_0-1)(v_1-1);\ l_0\not\in\Gamma \\ 0<l_1<(v_0-1)v_1-sv_0;\ l_1\not\in\Gamma \\
1\leq j\leq s\end{array}$} \\ \hline

\footnotesize{$\begin{array}{c}\Gamma=\langle 4,6,v_2\rangle\\ \mu=v_2+3 \\ \lambda_1=v_2-2\end{array}$} & \footnotesize{$v_2-2, v_2+2$} &  $v_2+1$ & $\begin{array}{c}-\frac{i}{12};\ i\in\{5,7,11,13\}\vspace{0.2cm}\\ -\frac{v_2+2(j-1)}{2v_2};\ 0\leq j<v_2;\ j\neq 1\end{array}$ \\ \hline
\end{tabular}
\vspace{3mm}
		
\noindent {\bf Table 1:} Analytical invariants for plane branches described in Proposition \ref{G=0}.
\end{center}

\section{ $\Sigma_{\Gamma,\Lambda}$-strata {\it versus} $\Sigma_{\Gamma,\mathcal{R}}$-strata} 

In what follows, we denote by $\Sigma_{\Gamma}$ the topological class of all plane branches with semigroup $\Gamma$. 

We set 
\begin{itemize}
\item $\Sigma_{\Gamma,\tau}$: the stratum of $\Sigma_{\Gamma}$ of all plane branches with Tjurina number $\tau$ fixed;
    \item $\Sigma_{\Gamma,\Lambda}$: the stratum of $\Sigma_{\Gamma}$ of all plane branches with set of values of $1$-forms $\Lambda$  fixed;
    \item $\Sigma_{\Gamma,\mathcal{R}}$: the stratum of $\Sigma_{\Gamma}$ of all plane branches with a set of roots of Bernstein polynomial given by $\mathcal{R}$. 
\end{itemize}

In \cite{He-99}, Hertling and Stahlke consider adjacency relations among the strata $\Sigma_{\Gamma,\tau}$ and  $\Sigma_{\Gamma,\mathcal{R}}$. Since $\tau=\mu-\sharp\Lambda\setminus\Gamma$, that is, $\Lambda$ is a finer analytic invariant than $\tau$, it is natural to compare these strata with $\Sigma_{\Gamma,\Lambda}$. 

These stratifications do not coincide and there are not inclusion relations among them.
We illustrate the relations between $\Sigma_{\Gamma,\Lambda}$-stratum and $\Sigma_{\Gamma,\mathcal{R}}$-stratum by taking the topological class given by the semigroup $\Gamma_f=\langle 6,7\rangle$ that is a case considered in
\cite{cassous-86}, \cite{cassous-87} and \cite{He-99}.

\begin{example}\label{strata}
Any plane branch $\mathcal{C}_f$ with value semigroup $\Gamma=\langle 6,7\rangle$ can be given by, according to \cite{zariski}, by
\[f=X_1^7+X_2^6+a_1X_1^5X_2^2+a_2X_1^4X_2^3+a_3X_2^3X_2^4+a_4X_1^5X_2^3+a_5X_1^4X_2^4+a_6X_1^5X_2^4.\]

Applying the algorithm presented in \cite{basestandard}, we compute the possible value set $\Lambda_f$. In the next table, we present the set $\mathcal{R}_0$ of all roots of $\widetilde{b}_f$ opposite to a spectral number, this is sufficient to determine the set $\mathcal{R}_f$ because the other $\mu_f-\sharp\mathcal{R}_0$ elements of $\mathcal{R}_f$ are given by $-(\sigma+1)$ with $
 \sigma\in\mathcal{S}_{f}\setminus\{-\rho;\ \rho\in\mathcal{R}_0\}$.

{\footnotesize
\begin{center}
\begin{tabular}{|c|c|c|c|}
\hline
Restrictions & $\Lambda_f\setminus\Gamma_f$ & $\tau_f$ &  $\mathcal{R}_0$ \\
\hline
$\begin{array}{c} a_1\neq 0\\ 9a_2^2-8a_1a_3+\frac{20}{7}a_1^3\neq 0\end{array}$ & $\{{\bf 15},{\bf 22},{\bf 23},{\bf 29}\}$ & $26$ & 
$\begin{array}{l}
\left \{ -\frac{\bf 15}{42}, -\frac{16}{42}, -\frac{17}{42}, -\frac{\bf 22}{42}, -\frac{\bf 23}{42}, -\frac{\bf 29}{42}  \right \} 
\ \mbox{if}\ a_2\neq 0, a_3\neq\frac{2a_1^2}{7} \vspace{0.2cm}\\

\left \{ -\frac{\bf 15}{42}, -\frac{16}{42}, -\frac{\bf 22}{42}, -\frac{\bf 23}{42}, -\frac{\bf 29}{42}  \right \} 
\ \mbox{if}\ a_2\neq 0, a_3=\frac{2a_1^2}{7} \vspace{0.2cm} \\

\left \{ -\frac{\bf 15}{42}, -\frac{17}{42}, -\frac{\bf 22}{42}, -\frac{\bf 23}{42}, -\frac{\bf 29}{42}  \right \} 
\ \mbox{if}\ a_2= 0, a_3\neq\frac{2a_1^2}{7}, a_3\neq\frac{5a_1^2}{14} \vspace{0.2cm} \\

\left \{ -\frac{\bf 15}{42}, -\frac{\bf 22}{42}, -\frac{\bf 23}{42}, -\frac{\bf 29}{42}  \right \} 
\ \mbox{if}\ a_2= 0, a_3=\frac{2a_1^2}{7}\vspace{0.1cm}

\end{array}$
\\
\hline

$\begin{array}{c} a_1= 0\\ a_2\neq 0\end{array}$ & $\{{\bf 16},{\bf 22},{\bf 23},{\bf 29}\}$ & $26$ & 
$\begin{array}{l}
\\ \left \{ -\frac{\bf 16}{42}, -\frac{17}{42}, -\frac{\bf 22}{42}, -\frac{\bf 23}{42}, -\frac{\bf 29}{42}  \right \} 
\ \mbox{if}\ a_3\neq 0 \vspace{0.2cm}\\

\left \{ -\frac{\bf 16}{42}, -\frac{\bf 22}{42}, -\frac{\bf 23}{42}, -\frac{\bf 29}{42}  \right \} 
\ \mbox{if}\ a_3=0\vspace{0.1cm}
\end{array}$
\\
\hline

$\begin{array}{c} a_1\neq 0\\ 9a_2^2-8a_1a_3+\frac{20}{7}a_1^3= 0\end{array}$ & $\{{\bf 15},{\bf 22},{\bf 29}\}$ & $27$ & 
$\begin{array}{l}
\left \{ -\frac{\bf 15}{42}, -\frac{16}{42}, -\frac{17}{42}, -\frac{\bf 22}{42}, -\frac{23}{42}, -\frac{\bf 29}{42}  \right \} 
\ \mbox{if}\ a_2\neq 0, a_3\neq\frac{2a_1^2}{7} \vspace{0.2cm}\\

\left \{ -\frac{\bf 15}{42}, -\frac{16}{42}, -\frac{\bf 22}{42}, -\frac{23}{42}, -\frac{\bf 29}{42}  \right \} 
\ \mbox{if}\ a_2\neq 0, a_3=\frac{2a_1^2}{7} \vspace{0.2cm} \\

\left \{ -\frac{\bf 15}{42}, -\frac{17}{42}, -\frac{\bf 22}{42}, -\frac{23}{42}, -\frac{\bf 29}{42}  \right \} 
\ \mbox{if}\ a_2= 0, a_3=\frac{5a_1^2}{14}, a_5\neq-\frac{5a_1^5}{2^4\cdot 3\cdot 7^3} \vspace{0.2cm} \\

\left \{ -\frac{\bf 15}{42}, -\frac{17}{42}, -\frac{\bf 22}{42}, -\frac{\bf 29}{42}  \right \} 
\ \mbox{if}\ a_2= 0, a_3=\frac{5a_1^2}{14}, a_5=-\frac{5a_1^5}{2^4\cdot 3\cdot 7^3}\vspace{0.1cm}

\end{array}$
\\
\hline

$\begin{array}{c} a_1=a_2= 0\\ a_3\neq 0\end{array}$ & $\{{\bf 17},{\bf 23},{\bf 29}\}$ & $27$ & 
$\begin{array}{l}
\\ 
\left \{  -\frac{\bf 17}{42}, -\frac{22}{42}, -\frac{\bf 23}{42}, -\frac{\bf 29}{42}  \right \} 
\ \mbox{if}\ a_4\neq 0 \vspace{0.2cm}\\

\left \{ -\frac{\bf 17}{42}, -\frac{\bf 23}{42}, -\frac{\bf 29}{42}  \right \} 
\ \mbox{if}\ a_4=0\vspace{0.1cm}
\end{array}$
\\
\hline

$\begin{array}{c} a_1=a_2=a_3= 0\\ a_4\neq 0\end{array}$ & $\{{\bf 22},{\bf 29}\}$ & $28$ & 

$\begin{array}{l}
\\
\left \{  -\frac{\bf 22}{42}, -\frac{23}{42}, -\frac{\bf 29}{42}  \right \} 
\ \mbox{if}\ a_5\neq 0 \vspace{0.2cm}\\

\left \{ -\frac{\bf 22}{42}, -\frac{\bf 29}{42}  \right \} 
\ \mbox{if}\ a_5=0\vspace{0.1cm}
\end{array}$
\\
\hline

$\begin{array}{c} a_1=a_2=a_3=a_4= 0\\ a_5\neq 0\end{array}$ & $\{{\bf 23},{\bf 29}\}$ & $28$ & 
$\left \{ -\frac{\bf 23}{42}, -\frac{\bf 29}{42}  \right \}$
\\
\hline

$\begin{array}{c} a_i= 0,\ 1\leq i\leq 5\\ a_6\neq 0\end{array}$ & $\{{\bf 29}\}$ & $29$ & 
$\left \{ -\frac{\bf 29}{42}  \right \}$
\\
\hline

$ a_i=0,\ 1\leq i\leq 6$ & $\emptyset$ & $30$ & 
$\emptyset$
\\
\hline
\end{tabular}

\end{center}
}

Notice that the last three lines of the above table satisfy the item 2 of Proposition \ref{G=0}, that is, the value set $\Lambda$ determines the set of roots of the Bernstein polynomial.
\end{example}

\vspace{0.5cm}
\noindent {\bf Acknowledgments}\\
The authors express their gratitude to Professor Abramo Hefez for his valuable discussions.

\vspace{0.75cm}
\begin{tabular}{lll}
{\sc Andrea Gomes Guimar\~ares} & \hspace{2cm} & {\sc Marcelo Escudeiro Hernandes} \\
Universidade Federal Fluminense & & Universidade Estadual de Maring\'a \\
Instituto de Matem\'atica e Estat\'istica & & Departamento de Matem\'atica \\
Niterói - RJ - Brazil & &   Maring\'a - PR - Brazil \\
angogui@gmail.com & & mehernandes@uem.br
\end{tabular}

\end{document}